\documentclass[generic, preprint]{imsart}

\setattribute{journal}{name}{}

\RequirePackage[OT1]{fontenc}
\RequirePackage{amsthm,amsmath,amsfonts}
\RequirePackage{graphicx}
\RequirePackage[numbers]{natbib}
\RequirePackage[colorlinks,citecolor=blue,urlcolor=blue]{hyperref}

\newtheorem{theorem}{Theorem}
\newtheorem{definition}{Definition}
\newtheorem{proposition}{Proposition}
\newtheorem{assumption}{Assumption}
\newtheorem{corollary}{Corollary}
\newtheorem{remark}{Remark}

\usepackage[margin=1in]{geometry}

\startlocaldefs
\newcommand{\bdsm}{\boldsymbol}
\endlocaldefs

\begin{document}

\begin{frontmatter}

\title{A Normal Hierarchical Model and Minimum Contrast Estimation for Random Intervals}
\runtitle{Normal Hierarchical MCE Random Intervals}



\begin{aug}
\author{\fnms{Yan} \snm{Sun}\corref{}\thanksref{}\ead[label=e1]{yan.sun@usu.edu}}
\and
\author{\fnms{Dan} \snm{Ralescu}\corref{}\thanksref{t2}\ead[label=e2]{dan.ralescu@uc.edu}}

\thankstext{t2}{Corresponding author. Supported by a Taft Research Grant.}

\runauthor{Y. Sun and D. Ralescu}

\affiliation{Utah State University and University of Cincinnati}

\address{Department of Mathematics $\&$ Statistics\\
Utah State University\\
3900 Old Main Hill\\
Logan, Utah 84322-3900\\
\printead{e1}\\
\phantom{E-mail: yan.sun@usu.edu}}

\address{Department of Mathematical Sciences\\
University of Cincinnati\\
4199 French Hall\\
2815 Commons Way\\
Cincinnati, Ohio 45221-0025\\
\printead{e2}\\
\phantom{E-mail: dan.ralescu@uc.edu}}
\end{aug}

\begin{abstract}
Many statistical data are imprecise due to factors such as measurement errors, computation errors, and lack of information. In such cases, data are better represented by intervals rather than by single numbers. Existing methods for analyzing interval-valued data include regressions in the metric space of intervals and symbolic data analysis, the latter being proposed in a more general setting. However, there has been a lack of literature on the parametric modeling and distribution-based inferences for interval-valued data.  In an attempt to fill this gap, we extend the concept of normality for random sets by Lyashenko and propose a Normal hierarchical model for random intervals. In addition, we develop a minimum contrast estimator (MCE) for the model parameters, which we show is both consistent and asymptotically normal. Simulation studies support our theoretical findings, and show very promising results. Finally, we successfully apply our model and MCE to a real dataset.
\end{abstract}



\end{frontmatter}


\section{Introduction}
In classical statistics, it is often assumed that the outcome of an experiment is precise and the uncertainty of observations is solely due to randomness. Under this assumption, numerical data are represented as collections of real numbers. In recent years, however, there has been increased interest in situations when exact outcomes of the experiment are very difficult or impossible to obtain, or to measure. The imprecise nature of the data thus collected is caused by various factors such as measurement errors, computational errors, loss or lack of information. Under such circumstances and, in general, any other circumstances such as grouping and censoring, when observations cannot be pinned down to single numbers, data are better represented by intervals.  Practical examples include interval-valued stock prices, oil prices, temperature data, medical records, mechanical measurements, among many others.

In the statistical literature, random intervals are most often studied in the framework of random sets, for which the probability-based theory has developed since the publication of the seminal book Matheron (1975). Studies on the corresponding statistical methods to analyze set-valued data, while still at the early stage, have shown promising advances. See Stoyan (1998) for a comprehensive review. Specifically, to analyze interval-valued data, the earliest attempt probably dates back to 1990, when Diamond published his paper on the least squares fitting of compact set-valued data and considered interval-valued input and output as a special case (see Diamond (1990)). Due to the embedding theorems started by Brunn and Minkowski and later refined by R\.adstr\"om (see R\.adstr\"om (1952)) and H\"ormander (see H\"ormander (1954)), $\mathcal{K}(\mathbb{R}^n)$, the space of all nonempty compact convex subsets of $\mathbb{R}^n$, is embedded into the Banach space of support functions. Diamond (1990) defined an $L_2$ metric in this Banach space of support functions, and found the regression coefficients by minimizing the $L_2$ metric of the sum of residuals. This idea was further studied in Gil et al. (2002), where the $L_2$ metric was replaced by a generalized metric on the space of nonempty compact intervals, called ``W-distance'', proposed earlier by K\"orner (1998). Separately, Billard and Diday (2003) introduced the central tendency and dispersion measures and developed the symbolic interval data analysis based on those. (See also Carvalho et al. (2004).) However, none of the existing literature considered distributions of the random intervals and the corresponding statistical methods.

It is well known that normality plays an important role in classical statistics. But the normal distribution for random sets remained undefined for a long time, until the 1980s when the concept of normality was first introduced for compact convex random sets in the Euclidean space by Lyashenko (1983). This concept is especially useful in deriving limit theorems for random sets. See, Puri et al. (1986), Norberg (1984), among others. Since a compact convex set in $\mathbb{R}$ is a closed bounded interval, by the definition of Lyashenko (1983), a normal random interval is simply a Gaussian displacement of a fixed closed bounded interval. From the point of view of statistics, this is not enough to fully capture the randomness of a general random interval.

In this paper, we extend the definition of normality given by Lyashenko (1983) and propose a Normal hierarchical model for random intervals. With one more degree of freedom on ``shape'', our model conveniently captures the entire randomness of random intervals via a few parameters. It is a natural extension from Lyashenko (1983) yet a highly practical model accommodating a large class of random intervals. In particular, when the length of the random interval reduces to zero, it becomes the usual normal random variable. Therefore, it can also be viewed as an extension of the classical normal distribution that accounts for the extra uncertainty added to the randomness. 
In addition, there are two interesting properties regarding our Normal hierarchical model: 1) conditioning on the first hierarchy, it is exactly the normal random interval defined by Lyashenko (1983), which could be a very useful property in view of the limit theorems; 2) with certain choices of the distributions, a linear combination of our Normal hierarchical random intervals follows the same Normal hierarchical distribution. An immediate consequence of the second property is the possibility of a factor model for multi-dimensional random intervals, as the ``factor'' will have the same distribution as the original intervals.

For random sets models, it is important, in the stage of parameter estimation, to take into account the geometric characteristics of the observations. For example, Tanaka et al. (2008) proposed an approximate maximum likelihood estimation for parameters in the Neyman-Scott point processes based on the point pattern of the observation window. For another model, Heinrich (1993) discussed several distance functions (called ``contrast functions'') between the parametric and the empirical contact distribution function that are used towards parameter estimation for Boolean models. Bearing this in mind, to estimate the parameters of our Normal hierarchical model, we propose a minimum contrast estimator (MCE) based on the hitting function (capacity functional) that characterizes the distribution of a random interval by the hit-and-miss events of test sets. See Matheron (1975). In particular, we construct a contrast function based on the integral of a discrepancy function between the empirical and the parametric distribution measure. Theoretically, we show that under certain conditions our MCE satisfies a strong consistency and asymptotic normality. The simulation study is consistent with our theorems. We apply our model to analyze a daily temperature range data and, in this context, we have derived interesting and promising results.

The use of an integral measure of probability discrepancy here is not new.  For example, the integral probability metrics (IPMs), widely used as tools for statistical inferences, have been defined as the supremum of the absolute differences between expectations with respect to two probability measures. See, e.g., Zolotarev (1983), M\"uller (1997), and Sriperumbudur et al. (2012), for references. Especially, the empirical estimation of IPMs proposed by Sriperumbudur et al. (2012) drastically reduces the computational burden, thereby emphasizing the practical use of the IPMs. This idea is potentially applicable to our MCE and we expect similar reduction in computational intensity as for IPMs.

The rest of the paper is organized as follows. Section \ref{sec:model} formally defines our Normal hierarchical model and discusses its statistical properties. Section \ref{sec:mce} introduces a minimum contrast estimator for the model parameters, and presents its asymptotic properties. A simulation study is reported in Section \ref{sec:simu}, and a real data application is demonstrated in Section \ref{sec:real}. We give concluding remarks in Section \ref{sec:conclu}. Proofs of the theorems are presented in Section \ref{sec:proofs}. Useful lemmas and other proofs are deferred to the Appendix.

\section{The Normal hierarchical model}\label{sec:model}
\subsection{Definition}
Let $(\Omega,\mathcal{L},P)$ be a probability space. Denote by
$\mathcal{K}$ the collection of all non-empty compact subsets of
$\mathbb{R}^d$. A random compact set is a Borel measurable function
$A: \Omega\rightarrow\mathcal{K}$, $\mathcal{K}$ being equipped with
the Borel $\sigma$-algebra induced by the Hausdorff metric. If
$A(\omega)$ is convex for almost all $\omega$, then $A$ is called a
random compact convex set. (See Molchanov (2005), p.21, p.102.)
Denote by $\mathcal{K}_{\mathcal{C}}$ the collection of all compact
convex subsets of $\mathbb{R}^d$. By Theorem 1 of Lyashenko (1983),
a compact convex random set $A$ in the Euclidean space
$\mathbb{R}^d$ is Gaussian if and only if $A$ can be represented as
the Minkowski sum of a fixed compact convex set $M$ and a
$d$-dimensional normal random vector $\epsilon$, i.e.
\begin{equation}\label{def_Lsko}
  A=M+\left\{\epsilon\right\}.
\end{equation}
As pointed out in Lyashenko (1983), Gaussian random sets are
especially useful in view of the limit theorems discussed earlier in
Lyashenko (1979). That is, if the conditions in those theorems are
satisfied and the limit exists, then it is Gaussian in the sense of
(\ref{def_Lsko}). Puri et al. (1986) extended these results to
separable Banach spaces.

In the following, we will restrict ourselves to compact convex random sets in $\mathbb{R}^1$, that is, bounded closed random intervals. They will be called random intervals for ease of presentation.

According to (\ref{def_Lsko}), a random interval $A$ is Gaussian if and only if A is representable in the form
\begin{equation}\label{def}
  A=I+\left\{\epsilon\right\},
\end{equation}
where $I$ is a fixed bounded closed interval and $\epsilon$ is a
normal random variable. Obviously, such a random interval is simply
a Gaussian displacement of a fixed interval, so it is not enough to
fully capture the randomness of a general random interval. In order
to model the randomness of both the location and the ``shape''
(length), we propose the following Normal hierarchical model for
random intervals:
\begin{align}
  A&= I+\left\{\epsilon\right\},\label{def:A_1}\\
  I&= \eta I_0,\label{def:A_2}
\end{align}
where $\eta$ is another random variable and $I_0=\left[a_0, b_0\right]$ is a fixed
interval in $\mathbb{R}$. Here, the product $\eta I_0$ is in the
sense of scalar multiplication of a real number and a set. Let
$\lambda(\cdot)$ denote the Lebesgue measure of $\mathbb{R}^1$.
 Then,
\begin{eqnarray}\label{length}
\lambda(A)=\lambda(\epsilon+\eta I_0)=\lambda(\eta I_0)=\left|\eta\right|\lambda(I_0).
\end{eqnarray}
That is, $\eta$ is the variable that models the length of $A$. In particular, if $\eta\to 0$, then A reduces to a normal random variable.

Obviously, $\epsilon$ and $\eta$ are ``location'' and ``shape''
variables. We assume that $\eta>0$. Then the Normal hierarchical
random interval is explicitly expressible as
\begin{equation*}
  A=\left[\epsilon+a_0\eta, \epsilon+b_0\eta\right].
\end{equation*}
The parameter $b_0$ is indeed unnecessary, as the difference $b_0-a_0$ can be
absorbed by $\eta$. As a result,
\begin{equation}\label{mod-simple}
  A=\left[\epsilon+a_0\eta, \epsilon+\left(a_0+1\right)\eta\right]
\end{equation}
Compared to the ``naive" model $A=[\epsilon-\frac{1}{2}\eta,
\epsilon+\frac{1}{2}\eta]$, for which $\epsilon$ is precisely the
center of the interval, (\ref{mod-simple}) has an extra parameter
$a_0$. Notice that the center of $A$ is
$\epsilon+\left(a_0+\frac{1}{2}\right)\eta$, so $a_0$ controls the
difference between $\epsilon$ and the center, and therefore is
interpreted as modeling the uncertainty that the Normal random variable $\epsilon$ is not
necessarily the center. 

\begin{remark}\label{rmk:1}
There are some existing works in the literature to model the randomness of intervals. For example, a random interval can be viewed as the ``crisp" version of the LR-fuzzy random variable, which is often used to model the randomness of imprecise intervals such as [approximately 2, approximately 5]. See K\"orner (1997) for detailed descriptions. However, as far as the authors are aware, models with distribution assumptions for interval-valued data have not been studied yet. Our Normal hierarchical random interval is the first statistical approach that extends the concept of normality while modeling the full randomness of an interval.
\end{remark}

An interesting property of the Normal hierarchical random interval is that its linear combination is still a Normal hierarchical random interval. This is seen by simply observing that
\begin{eqnarray}\label{equ:factor}
\sum\limits_{i=1}^{n}a_iA_i=\sum\limits_{i=1}^{n}a_i\left(\epsilon_i+\eta_iI_0\right)
=\sum\limits_{i=1}^{n}a_i\epsilon_i+I_0\left(\sum\limits_{i=1}^{n}a_i\eta_i\right),
\end{eqnarray}
for arbitrary constants $a_i, i=1,\cdots,n$, where ``$+$'' denotes the Minkowski addition. This is very useful in developing a factor model for the analysis of multiple random intervals. Especially, if we assume $\eta_i\sim N(\mu_i,\sigma^2_i), i=1,\cdots,n$, then the ``factor'' $\sum\limits_{i=1}^{n}a_iA_i$ has exactly the same distribution as the original random intervals. We will elaborate more on this issue in section \ref{sec:simu}.

Without loss of generality, we can assume in the model (\ref{def:A_1})-(\ref{def:A_2}) that $E\epsilon=0$. We will make this assumption throughout the rest of the paper.

\subsection{Model properties}
According to the Choquet theorem ( Molchanov (2005), p.10), the distribution of a random closed set 
(and random compact convex set as a special case) A, 
is completely characterized by the
hitting function $T$ defined as:
\begin{equation}\label{dist}
  T(K)=P(K\cap A\neq\emptyset), \ \ \forall K\in\mathcal{K}_{\mathcal{C}}.
\end{equation}

Writing $I_0=[a_0,b_0]$ with $a_0\leq b_0$, the Normal hierarchical
random interval in (\ref{def:A_1})-(\ref{def:A_2}) has the following
hitting function: for $K=[a,b]$:
\begin{eqnarray*}
  &&T_A([a,b])\\
  &=&P([a,b]\cap A\neq\emptyset)\\
  &=&P([a,b]\cap A\neq\emptyset,\eta\geq 0)+P([a,b]\cap A\neq\emptyset,\eta< 0)\\
  &=&P(a-\eta b_0\leq\epsilon\leq b-\eta a_0,\eta\geq 0)+P(a-\eta a_0\leq\epsilon\leq b-\eta b_0,\eta< 0).
\end{eqnarray*}

The expectation of a compact convex random set $A$ is defined by the
Aumann integral (see Aumann (1965), Artstein and Vitale (1975)) as
\begin{eqnarray*}
  EA=\left\{E\xi:\xi\in A \text{ almost surely}\right\}.
\end{eqnarray*}
In particular, the Aumann expectation of a random interval $A$ is given by
\begin{equation}\label{def:aumann}
  EA=[EA_l,EA_u],
\end{equation}
where $A_l$ and $A_u$ are the interval ends. Therefore, the Aumann
expectation of the Normal hierarchical random interval $A$ is
\begin{eqnarray*}
  EA&=&E(\epsilon+\eta I_0)=E\epsilon+E(\eta I_0)=E(\eta I_0)\\
    &=&E\left\{[a_0\eta,b_0\eta]I_{(\eta\geq 0)}+[b_0\eta,a_0\eta]I_{(\eta<0)}\right\}\\
    &=&E\left[a_0\eta I_{(\eta\geq 0)}+b_0\eta I_{(\eta<0)},b_0\eta I_{(\eta\geq 0)}+a_0\eta I_{(\eta<0)}\right]\\
    &=&\left[a_0E\eta_{+}+b_0E\eta_{-},b_0E\eta_{+}+a_0E\eta_{-}\right],
\end{eqnarray*}
where
\begin{eqnarray*}
  \eta_{+}&=&\eta I_{(\eta\geq 0)},\\
  \eta_{-}&=&\eta I_{(\eta< 0)}.
\end{eqnarray*}
Notice that $\eta_{+}$ can be interpreted as the positive part of $\eta$, but $\eta_{-}$ is not the negative part of $\eta$, as $\eta_{-}<0$ when $\eta<0$.

The variance of a compact convex random set $A$ in $\mathbb{R}^d$ is
defined via its support function. 
In the special case when $d=1$, it is shown by
straightforward calculations that
\begin{equation}\label{var-1}
  Var(A)=\frac{1}{2}Var(A_l)+\frac{1}{2}Var(A_u),
\end{equation}
 or equivalently,
\begin{equation}\label{var-2}
  Var(A)=Var\left(A_c\right)+Var\left(A_r\right),
\end{equation}
where $A_c$ and $A_r$ denote the center and radius of a random interval $A$. See K\"orner (1995). Again, as we pointed out in Remark \ref{rmk:1}, a random interval can be viewed as a special case of the LR-fuzzy random variable. Therefore, formulae (\ref{var-1}) and (\ref{var-2}) coincide with the variance of the LR-fuzzy random variable, when letting the left and right spread both equal to 0, i.e., $l=r=0$. See K\"orner (1997). For the Normal hierarchical random interval $A$,
\begin{eqnarray*}
  &&Var(A_l)\\
  &=&Var\left(\epsilon+a_0\eta_{+}+b_0\eta_{-}\right)\\
  &=&E\left(\epsilon+a_0\eta_{+}+b_0\eta_{-}\right)^2-\left[E\left(\epsilon+a_0\eta_{+}+b_0\eta_{-}\right)\right]^2\\
  &=&E\epsilon^2+a_0^2Var(\eta_{+})+b_0^2Var(\eta_{-})\\
  &&+2\left(a_0E\epsilon\eta_{+}+b_0E\epsilon\eta_{-}-a_0b_0E\eta_{+}E\eta_{-}\right),
\end{eqnarray*}
and, analogously,
\begin{eqnarray*}
  &&Var(A_u)\\
  &=&E\epsilon^2+b_0^2Var(\eta_{+})+a_0^2Var(\eta_{-})\\
  &&+2\left(b_0E\epsilon\eta_{+}+a_0E\epsilon\eta_{-}-a_0b_0E\eta_{+}E\eta_{-}\right).
\end{eqnarray*}
The variance of $A$ is then found to be
\begin{eqnarray*}
  Var(A)&=&\frac{1}{2}Var(A_l)+\frac{1}{2}Var(A_u)\\
  &=&E\epsilon^2+\frac{1}{2}\left(a_0^2+b_0^2\right)\left[Var(\eta_{+})+Var(\eta_{-})\right]\\
  &&+(a_0+b_0)E\epsilon\eta-2a_0b_0E\eta_{+}\eta_{-}.
\end{eqnarray*}

\begin{remark}
Assuming $\eta>0$, we have
\begin{eqnarray*}
  Var(A)
  &=&E\epsilon^2+\frac{1}{2}(a_0^2+b_0^2)Var(\eta)+(a_0+b_0)E\epsilon\eta\\
  &=&Var(\epsilon)+\frac{1}{2}(a_0^2+b_0^2)Var(\eta)+(a_0+b_0)Cov(\epsilon,\eta),
\end{eqnarray*}
with $E\epsilon=0$. This formula certainly includes the special case
of the ``naive" model $A=[\epsilon-\frac{1}{2}\eta,
\epsilon+\frac{1}{2}\eta]$, by letting $a_0=-\frac{1}{2}$ and
$b_0=\frac{1}{2}$. It is more general because it also accounts for
the covariance between ``location" and ``length" in calculating the
total variance of the random interval, while the ``naive" model
simply has
$Var\left(A\right)=Var\left(\epsilon\right)+Var\left(\eta\right)$.
\end{remark}

\section{The minimum contrast estimation}\label{sec:mce}
\subsection{Definitions}
We study minimum contrast estimation (MCE) of the parameters of the
Normal hierarchical random interval
($\ref{def:A_1}$)-($\ref{def:A_2}$), as well as its asymptotic
properties. Since $d=1$, from now on we let $\mathcal{K}$ be the
space of all non-empty compact subsets in $\mathbb{R}$
restrictively, and let $\mathcal{F}$ be the Borel $\sigma$-algebra
on $\mathcal{K}$ induced by the Hausdorff metric. Let
$\mathcal{K}_{\mathcal{C}}$ denote the space of all non-empty
compact convex subsets, i.e., bounded closed intervals, in
$\mathbb{R}$. As mentioned in the previous section, a random
interval $X$ is a Borel measurable function from a probability space
$(\Omega,\mathcal{L},P)$ to $(\mathcal{K},\mathcal{F})$ such that
$X\in\mathcal{K}_{\mathcal{C}}$ almost surely.

Throughout this section, we assume observing a sample of i.i.d.
random intervals $X(n)=\left\{X_1,X_2,\cdots,X_n\right\}$. Let
$\bdsm{\theta}$ denote a $p\times 1$ vector containing all the
parameters in the model, which takes on a value from a parameter
space $\Theta\subset\mathbb{R}^p$. Here $p$ is the number of
parameters. Let $\bdsm{\theta}_0$ denote the true value of the
parameter vector. Denote by $T_{\bdsm{\theta}}([a,b])$
the hitting function of $X_i$ with parameter $\bdsm{\theta}$.

In order to introduce the MCE, we will need some extra notations. Let $\textbf{X}$ be a basic set and $\mathcal{A}$ be a $\sigma$-field over it. Let $\mathcal{B}$ denote a family of probability measures on (\textbf{X},$\mathcal{A}$) and $\tau$ be a mapping from $\mathcal{B}$ to some topologial space $T$. $\tau(P)$ denotes the parameter value pertaining to $P$, $\forall P\in\mathcal{B}$. The classical definition of MCE given in Pfanzagl (1969) is quoted below.

\begin{definition}
$\left[Pfanzagl (1969)\right]$ A family of $\mathcal{A}$-measurable functions $f_t:\textbf{X}\rightarrow\mathbb{R},t\in T$ is a family of contrast functions if
\begin{equation}\label{cond1}
  E_P\left[f_t\right]<\infty,
\end{equation}
$\forall t\in T, \forall P\in\mathcal{B}$, and
\begin{equation}\label{cond2}
  E_P\left[f_{\tau(P)}\right]<E_P\left[f_t\right],
\end{equation}
$\forall t\in T, \forall P\in\mathcal{B}, t\neq\tau(P)$.
\end{definition}

In other words, a contrast function is a measurable function of the
random variable(s) whose expected value reaches its minimum under
the probability measure that generates the random variable(s). From
the view of probability, with the true parameters, a contrast
function tends to have a smaller value than with other parameters.

Adopting notation from Pfanzagl (1969), we let $\mathcal{B}$ denote
a family of probability measures on
($\mathcal{K}_{\mathcal{C}},\mathcal{F}$) and $\tau$ be a mapping
from $\mathcal{B}$ to some topologial space $T$. Similarly,
$\tau(P)$ denotes the parameter value pertaining to $P$, $\forall
P\in\mathcal{B}$. In a similar fashion to the contrast
function in Heinrich (1993) for Boolean models, we give our
definition of contrast function for random intervals in the
following. And then the MCE is defined as the minimizer of the
contrast function.

\begin{definition}\label{def:cf}
A family of $\mathcal{F}^n$-measurable functions $M(X(n);\boldsymbol{\theta})$: $\mathcal{K}_{\mathcal{C}}^n\rightarrow [-\infty,+\infty]$, $n\in\mathbb{N}$, $\boldsymbol{\theta}\in\boldsymbol{\Theta}$ is a family of contrast functions for $\mathcal{B}$, if there exists a function $N(\cdot,\cdot)$: $\boldsymbol{\Theta}\times\boldsymbol{\Theta}\rightarrow\mathbb{R}$ such that
\begin{equation}
  P_{\boldsymbol{\theta}}(\left\{\omega: \lim_{n\rightarrow\infty}M(X(n);\boldsymbol{\zeta})
  =N(\boldsymbol{\theta},\boldsymbol{\zeta})\right\})=1,\ \forall\ \boldsymbol{\theta},\boldsymbol{\zeta}
  \in\boldsymbol{\Theta},
\end{equation}
and
\begin{equation}
  N(\boldsymbol{\theta}, \boldsymbol{\theta})<N(\boldsymbol{\theta}, \boldsymbol{\zeta})\ \forall\
  \boldsymbol{\theta}, \boldsymbol{\zeta}\in\boldsymbol{\Theta},\ \boldsymbol{\theta}\neq\boldsymbol{\zeta}.
\end{equation}
\end{definition}

\begin{definition}\label{def:mce}
A $\mathcal{F}^n$-measurable function $\hat{\boldsymbol{\theta}}_n$: $\mathcal{K}_{\mathcal{C}}^n\rightarrow
\tau(\mathcal{B})$, which depends on $X(n)$ only, is called a minimum contrast estimator (MCE) if
\begin{equation}
  M(X(n);\hat{\boldsymbol{\theta}}_n)=\inf\left\{M(X(n);\boldsymbol{\theta}): \boldsymbol{\theta}\in
  \tau(\mathcal{B})\right\}.
\end{equation}
\end{definition}

\subsection{Theoretical results}
We make the following assumptions to present the theoretical results in this section.
\begin{assumption}\label{aspt:1}
$\Theta$ is compact, and $\bdsm{\theta}_0$ is an interior point of $\Theta$.
\end{assumption}
\begin{assumption}\label{aspt:2}
The model is identifiable.
\end{assumption}
\begin{assumption}\label{aspt:3}
$T_{\bdsm{\theta}}([\cdot,\cdot])$ is continuous with respect to $\bdsm{\theta}$.
\end{assumption}
\begin{assumption}\label{aspt:4}
$\frac{\partial T_{\bdsm{\theta}_0}}{\partial \theta_i}([\cdot,\cdot])$, $i=1,\cdots,p$, exist and are finite on a bounded region $S^0\subset\mathbb{R}^2$.
\end{assumption}
\begin{assumption}\label{aspt:5}
$\frac{\partial T_{\bdsm{\theta}}}{\partial\theta_j}([\cdot,\cdot])$, $\frac{\partial^2T_{\bdsm{\theta}}}{\partial\theta_j\partial\theta_k}([\cdot,\cdot])$, and
$\frac{\partial^3T_{\bdsm{\theta}}}{\partial\theta_j\partial\theta_k\partial\theta_l}([\cdot,\cdot])$, $i,j,k=1,\cdots,p$, exist and are finite on $S^0$ for $\bdsm{\theta}\in\Theta$.
\end{assumption}

Assumptions 4 and 5 are essential to establish the asymptotic normality for the MCE $\hat{\boldsymbol{\theta}}_n$. They are rather mild and can be met by a large class of capacity functionals. For example, if $S^0$ is closed, then each $T_{\bdsm{\theta}_0}$ with continuous up to third order partial derivatives satisfies both assumptions, as a continuous function on a compact region is always bounded. The following theorem gives sufficient conditions under which the minumum contrast estimator $\hat{\boldsymbol{\theta}}_n$ defined above is strongly consistent.
\begin{theorem}\label{thm:strong-consist}
Let $M(X(n);\boldsymbol{\theta})$ be a contrast function as in Definition \ref{def:cf} and let $\hat{\boldsymbol{\theta}}_n$ be the corresponding MCE. Under the hypothesis of Assumption \ref{aspt:1} and in addition if $M(X(n);\boldsymbol{\theta})$ is equicontinuous w.r.t. $\boldsymbol{\theta}$ for all $X(n), n=1,2,\cdots$, then,
\begin{equation*}
  \hat{\boldsymbol{\theta}}_n\rightarrow \boldsymbol{\theta}_0\ \ a.s.,\ \text{as}\  n\rightarrow\infty.
\end{equation*}
\end{theorem}

Let $[a,b]\in\mathcal{K}_{\mathcal{C}}$. Define an empirical estimator $\hat{T}([a,b];X(n))$ for $T([a,b])$ as:
\begin{equation}\label{T_hat}
  \hat{T}([a,b];X(n))=\frac{\# \left\{X_i: [a,b]\cap X_i\neq\emptyset, i=1,\cdots,n\right\}}{n}.
\end{equation}
Extending the contrast function defined in Heinrich (1993) (for parameters in the Boolean model), we construct a family of functions:
\begin{equation}\label{H_def}
  H(X(n);\boldsymbol{\theta})
  =\iint\limits_{S}\left[T_{\boldsymbol{\theta}}([a,b])-\hat{T}([a,b];X(n))\right]^2W(a,b)\mathrm{d}a\mathrm{d}b,
\end{equation}
for $\boldsymbol{\theta}\in\boldsymbol{\Theta}$, where $S\subset S^0\subset\mathbb{R}^2$, and $W(a,b)$ is a weight function on $[a,b]$ satisfying $0<W(a,b)<C$, $\forall [a,b]\in\mathcal{K}_\mathcal{C}$.

We show in the next Proposition that $H(X(n);\bdsm{\theta})$, $\bdsm{\theta}\in\bdsm{\Theta}$ defined in (\ref{H_def}) is a family of contrast functions for $\bdsm{\theta}$. This, together with Theorem \ref{thm:strong-consist}, immediately yields the strong consistency of the associated MCE. This result is summarized in Corollary \ref{coro:consist}.
\begin{proposition}\label{prop:cf}
Suppose that Assumption \ref{aspt:2} and Assumption \ref{aspt:3} are satisfied. Then $H(X(n);\boldsymbol{\theta})$, $\boldsymbol{\theta}\in\boldsymbol{\Theta}$, as defined in (\ref{H_def}), is a family of contrast functions with limiting function
\begin{equation}
  N(\boldsymbol{\theta},\boldsymbol{\zeta})=
  \iint\limits_{S}\left[T_{\boldsymbol{\theta}}([a,b])-T_{\boldsymbol{\zeta}}([a,b])\right]^2W(a,b)\mathrm{d}a\mathrm{d}b.
\end{equation}
In addition, $H(X(n);\boldsymbol{\theta})$ is equicontinuous w.r.t. $\boldsymbol{\theta}$.
\end{proposition}

\begin{corollary}\label{coro:consist}
Suppose that Assumption \ref{aspt:1}, Assumption \ref{aspt:2}, and Assumption \ref{aspt:3} are satisfied. Let $H(X(n);\boldsymbol{\theta})$ be defined as in (\ref{H_def}), and
\begin{equation}\label{def:theta}
  \bdsm{\theta}_n^H=\arg\min_{\bdsm{\theta}\in\Theta}H\left(X(n);\bdsm{\theta}\right).
\end{equation}
Then
\begin{equation}
  \bdsm{\theta}_n^H\rightarrow\bdsm{\theta}_0,\ a.s.,\nonumber
\end{equation}
as $n\rightarrow\infty$.
\end{corollary}

Next, we show the asymptotic normality for $\bdsm{\theta}_n^H$. As a preparation, we first prove the following proposition. The central limit theorem for $\bdsm{\theta}_n^H$ is then presented afterwards.
\begin{proposition}\label{prop:parH}
Assume the conditions of Lemma 1 (in the Appendix). Define
\begin{equation}
 \frac{\partial H}{\partial\bdsm{\theta}}\left(X(n);\bdsm{\theta}\right)
 :=\left[\frac{\partial H}{\partial\theta_1}\left(X(n);\bdsm{\theta}\right),\cdots,
 \frac{\partial H}{\partial\theta_p}\left(X(n);\bdsm{\theta}\right)\right]^{T},\nonumber
\end{equation}
as the $p\times 1$ gradient vector of $H\left(X(n);\bdsm{\theta}\right)$ w.r.t. $\bdsm{\theta}$. Then,
\begin{equation}
 \sqrt{n}\left[\frac{\partial H}{\partial\bdsm{\theta}}\left(X(n);\bdsm{\theta}_0\right)\right]
 \stackrel{\mathcal{D}}{\rightarrow}N\left(0,\Xi\right),\nonumber
\end{equation}
where $\Xi$ is the $p\times p$ symmetric matrix with the $(i,j)^{\text{th}}$ component
\begin{eqnarray}
 \Xi(i,j)&=&4\iiiint\limits_{S\times S}\left\{P\left(X_1\cap[a,b]\neq\emptyset,X_1\cap[c,d]\neq\emptyset\right)
  -T_{\bdsm{\theta}_0}\left([a,b]\right)T_{\bdsm{\theta}_0}\left([c,d]\right)\right\}\nonumber\\
  &&\frac{\partial T_{\bdsm{\theta}_0}}{\partial\theta_i}\left([a,b]\right)
  \frac{\partial T_{\bdsm{\theta}_0}}{\partial\theta_j}\left([c,d]\right)
  W(a,b)W(c,d)\mathrm{d}a\mathrm{d}b\mathrm{d}c\mathrm{d}d.\label{def:xi}
\end{eqnarray}
\end{proposition}

\begin{theorem}\label{thm:clt}
Let $H(X(n);\bdsm{\theta})$ be defined in (\ref{H_def}) and $\bdsm{\theta}_n^H$ be defined in (\ref{def:theta}). Assume the conditions of Corollary \ref{coro:consist}. If additionally Assumption \ref{aspt:5} is satisfied, then
\begin{equation}
 \sqrt{n}\left(\bdsm{\theta}_n^H-\bdsm{\theta}_0\right)\stackrel{\mathcal{D}}{\rightarrow}
 N\left(0,C(T_{\bdsm{\theta}_0})^{-1}\Xi C(T_{\bdsm{\theta}_0})^{-1}\right),
\end{equation}
where $C(T_{\bdsm{\theta}_0})=2\iint\limits_{S}\left(\frac{\partial T_{\bdsm{\theta}_0}}{\partial\bdsm{\theta}}\right)
 \left(\frac{\partial T_{\bdsm{\theta}_0}}{\partial\bdsm{\theta}}\right)^{T}([a,b])W(a,b)\mathrm{d}a\mathrm{d}b$, and $\Xi$ is defined in (\ref{def:xi}).
\end{theorem}

\section{Simulation}\label{sec:simu}
We carry out a small simulation to investigate the performance of the MCE introduced in Definition \ref{def:mce}. Assume, in the Normal hierarchical model (\ref{def:A_1})-(\ref{def:A_2}), that
\begin{equation}
  \begin{bmatrix} \epsilon\\ \eta \end{bmatrix}\sim
  \text{BVN}\left(
  \begin{bmatrix}0\\ \mu \end{bmatrix},
  \Sigma=\begin{bmatrix}\sigma_1^2 & \sigma_{12}\\ \sigma_{12} & \sigma_2^2 \end{bmatrix}
  \right),
 \end{equation}
  and
\begin{equation}
  b_0=a_0+1.
\end{equation}
The bivariate normal distribution conveniently takes care of the variances and covariance of the location variable $\epsilon$ and the shape variable $\eta$. 
The removal of the freedom of $b_0$ is for model identifiability purposes; it is seen that the hitting function $T_A$ is defined via $\eta a_0$ and $\eta b_0$ only. For the simulation, we assign the following parameter values:
\begin{equation}\label{eqn:par-val}
  a_0=1, \mu=20, \Sigma=\begin{bmatrix}10 & 1\\ 1 & 10\end{bmatrix}.
\end{equation}

\subsection{Hitting function}
Under the bivariate normal distribution assumption, the hitting function of our Normal hierarchical model is found to be
\begin{eqnarray}
  &&T_{\bdsm{\theta}}([a,b])\nonumber\\
  &=&P(a-\eta b_0\leq\epsilon\leq b-\eta a_0,\eta\geq 0)+P(a-\eta a_0\leq\epsilon\leq b-\eta b_0,\eta< 0)\nonumber\\
  &=&P\left(\epsilon\leq b-\eta a_0,\eta\geq 0\right)-P\left(\epsilon< a-\eta b_0,\eta\geq 0\right)\nonumber\\
  &&+P\left(\epsilon\leq b-\eta b_0,\eta< 0\right)-P\left(\epsilon< a-\eta a_0,\eta<0\right)\nonumber\\
  &=&P\left(\begin{bmatrix}1 & a_0\\ 0 & -1\end{bmatrix} \begin{bmatrix}\epsilon\\ \eta\end{bmatrix}\leq\begin{bmatrix}b\\0\end{bmatrix}\right)
  -P\left(\begin{bmatrix}1 & b_0\\ 0 & -1\end{bmatrix} \begin{bmatrix}\epsilon\\ \eta\end{bmatrix}\leq\begin{bmatrix}a\\0\end{bmatrix}\right)\nonumber\\
  &&+P\left(\begin{bmatrix}1 & b_0\\ 0 & 1\end{bmatrix} \begin{bmatrix}\epsilon\\ \eta\end{bmatrix}\leq\begin{bmatrix}b\\0\end{bmatrix}\right)
  -P\left(\begin{bmatrix}1 & a_0\\ 0 & 1\end{bmatrix} \begin{bmatrix}\epsilon\\ \eta\end{bmatrix}\leq\begin{bmatrix}a\\0\end{bmatrix}\right)\nonumber\\
  &=&\Phi\left(\begin{bmatrix}b\\0\end{bmatrix}; D_1\begin{bmatrix}0\\ \mu\end{bmatrix}, D_1\Sigma D_1^{'}\right)
  -\Phi\left(\begin{bmatrix}a\\0\end{bmatrix}; D_2\begin{bmatrix}0\\ \mu\end{bmatrix}, D_2\Sigma D_2^{'}\right)\nonumber\\
  &&+\Phi\left(\begin{bmatrix}b\\0\end{bmatrix}; D_3\begin{bmatrix}0\\ \mu\end{bmatrix}, D_3\Sigma D_3^{'}\right)
  -\Phi\left(\begin{bmatrix}a\\0\end{bmatrix}; D_4\begin{bmatrix}0\\ \mu\end{bmatrix}, D_4\Sigma D_4^{'}\right),\label{eqn:hit-fct}
\end{eqnarray}
where $\Phi\left(\textbf{x}; \bdsm{\mu}, \Omega\right)$ is the bivariate normal cdf with mean $\bdsm{\mu}$ and covariance $\Omega$, and
\begin{eqnarray*}
  D_1=\begin{bmatrix}1 & a_0\\ 0 & -1\end{bmatrix}, D_2=\begin{bmatrix}1 & b_0\\ 0 & -1\end{bmatrix},
  D_3=\begin{bmatrix}1 & b_0\\ 0 & 1\end{bmatrix}, D_4=\begin{bmatrix}1 & a_0\\ 0 & 1\end{bmatrix}.
\end{eqnarray*}
After linear transformation of variables, the terms in formula
(\ref{eqn:hit-fct}) is calculated via the standard bivariate normal
cdf. By absolute continuity, $T_{\bdsm{\theta}}([a,b])$ in this case
is continuous and also infinitely continuously differentiable.
Therefore, all the assumptions are satisfied and the corresponding
MCE achieves the strong consistency and asymptotic normality.

According to the assigned parameter values given in (\ref{eqn:par-val}), $P(\eta<0)<10^{-10}$. Therefore the hitting function is well approximated by
\begin{eqnarray*}
  &&T_{\bdsm{\theta}}([a,b])\\
  &\approx&P(a-\eta b_0\leq\epsilon\leq b-\eta a_0,\eta\geq 0)\\
  &\approx&P(a-\eta b_0\leq\epsilon\leq b-\eta a_0)\\
  &=&P\left(
  \begin{bmatrix}1 & a_0\\ -1 & -a_0-1\end{bmatrix}
  \begin{bmatrix}\epsilon\\ \eta\end{bmatrix}\leq
  \begin{bmatrix}b\\-a\end{bmatrix}\right)\\
  &=&\Phi\left(
  \begin{bmatrix}b\\-a\end{bmatrix};
  D\begin{bmatrix}0\\ \mu\end{bmatrix}, D\Sigma D^{'}\right),
\end{eqnarray*}
where
\begin{equation*}
  D=\begin{bmatrix}1 & a_0\\ -1 & -a_0-1\end{bmatrix}.
\end{equation*}
We use this approximate hitting function to simplify computation in our simulation study.


\subsection{Parameter initialization}
The model parameters can be estimated by the method of
moments. In most cases it is reasonable to assume $\eta^{-}\approx
0$, and consequently, $\eta\approx|\eta|$. So the moment estimates
for $\mu$ and $a_0$ are approximately
\begin{eqnarray}
 && \tilde{\mu} \leftarrow \bar{X_u}-\bar{X_l},\label{mm-1}\\
 && \tilde{a_0} \leftarrow \bar{X_l}/\tilde{\mu},\label{mm-2}
\end{eqnarray}
where $\bar{X_u}$ and $\bar{X_l}$ denote the sample means of $A_u$ and $A_l$, respectively. Denoting by $A_c$ the center of the random interval $A$, we further notice that $A_c=\epsilon+\frac{1}{2}\left(a_0+b_0\right)\eta=\epsilon+\left(a_0+\frac{1}{2}\right)\eta$. By the same approximation we have $\epsilon\approx A_c-\left(a_0+\frac{1}{2}\right)|\eta|$. Define a random variable
\begin{equation*}
  A_\delta=A_c-\left(a_0+\frac{1}{2}\right)|\eta|.
\end{equation*}
Then, the moment estimate for $\Sigma$ is approximately given by the sample variance-covariance matrix of $A_\delta$ and $A_u-A_l$, i.e.
\begin{equation}\label{mm-3}
  \tilde{\Sigma} \leftarrow \Sigma_s\left(A_\delta, A_u-A_l\right).
\end{equation}

\subsection{Performance of MCE}

Our simulation experiment is designed as follows: we first simulate an i.i.d. random sample of size $n$ from model (\ref{def:A_1})-(\ref{def:A_2}) with the assigned parameter values, then find the initial parameter values by (\ref{mm-1})-(\ref{mm-3}) based on the simulated sample, and lastly the initial values are updated to the MCE using the function \textsl{fminsearch.m} in Matlab 2011a. The process is repeated 10 times independently for each $n$, and we let $n=100, 200, 300, 400, 500$, successively, to study the consistency and efficiency of the MCE's.

\begin{figure}[ht]
\centering
\includegraphics[ height=2.000in, width=2.500in]{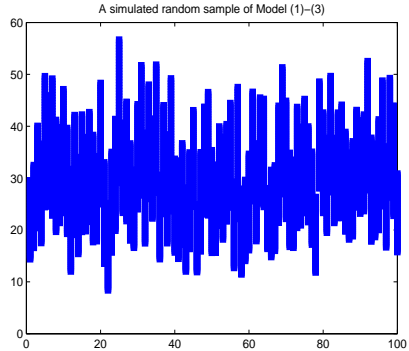}
\caption{Plot of a simulated sample from model (\ref{def:A_1})-(\ref{def:A_2}) with $n=100$.}
\label{fig:sample_simu}
\end{figure}

\begin{figure}[ht]
\centering
\includegraphics[ height=1.500in, width=2.300in]{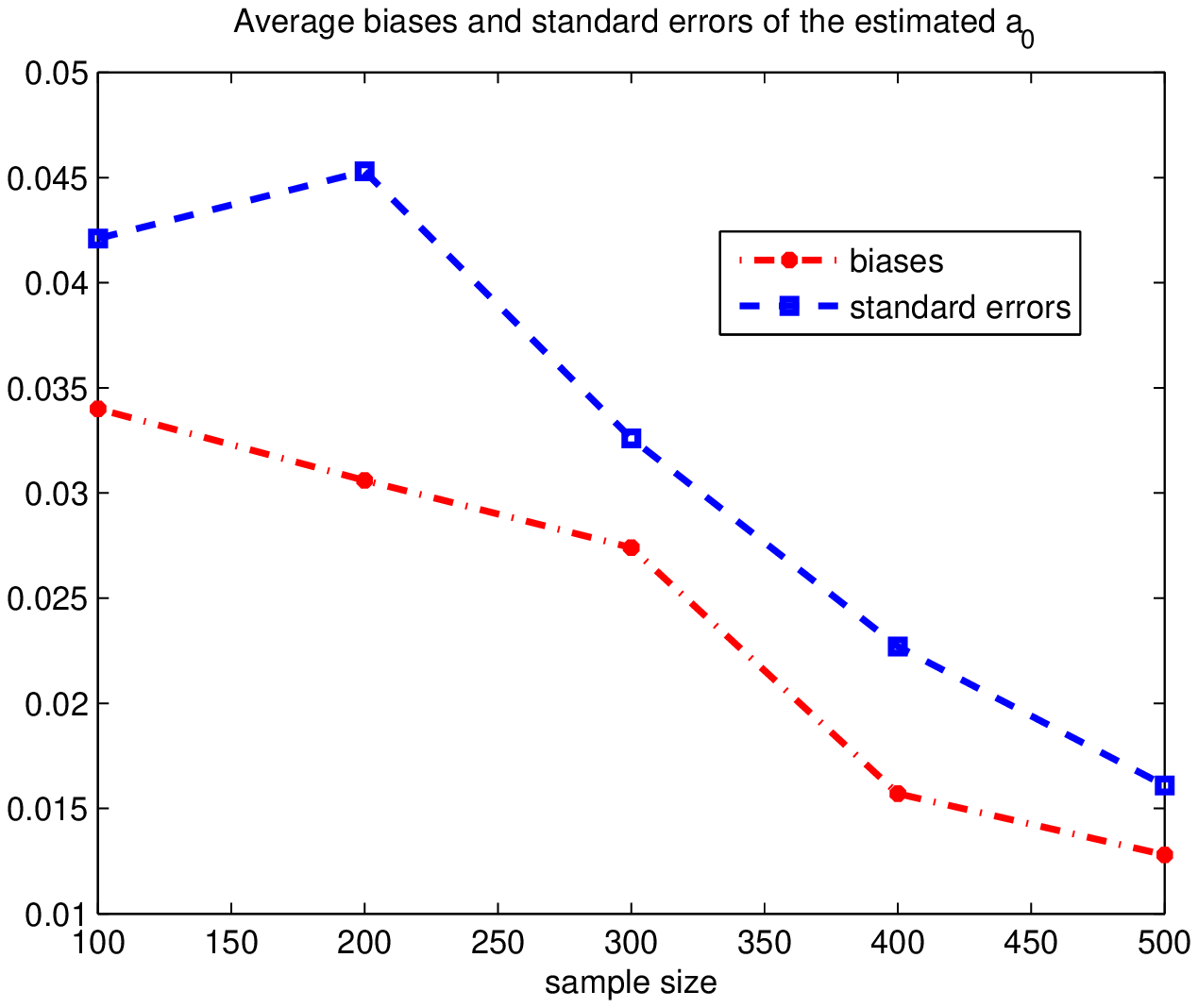}
\includegraphics[ height=1.500in, width=2.300in]{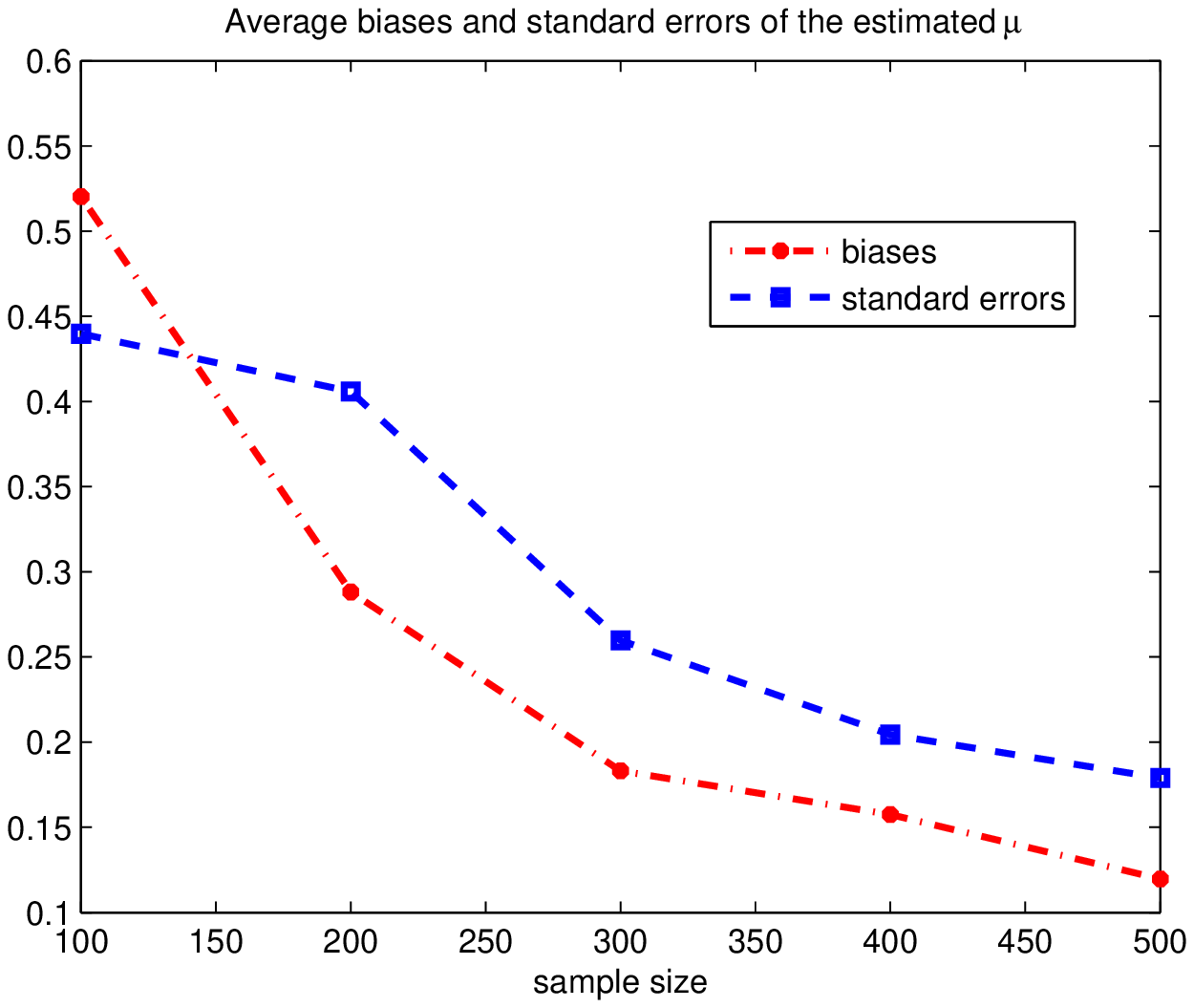}\\
\includegraphics[ height=1.500in, width=2.300in]{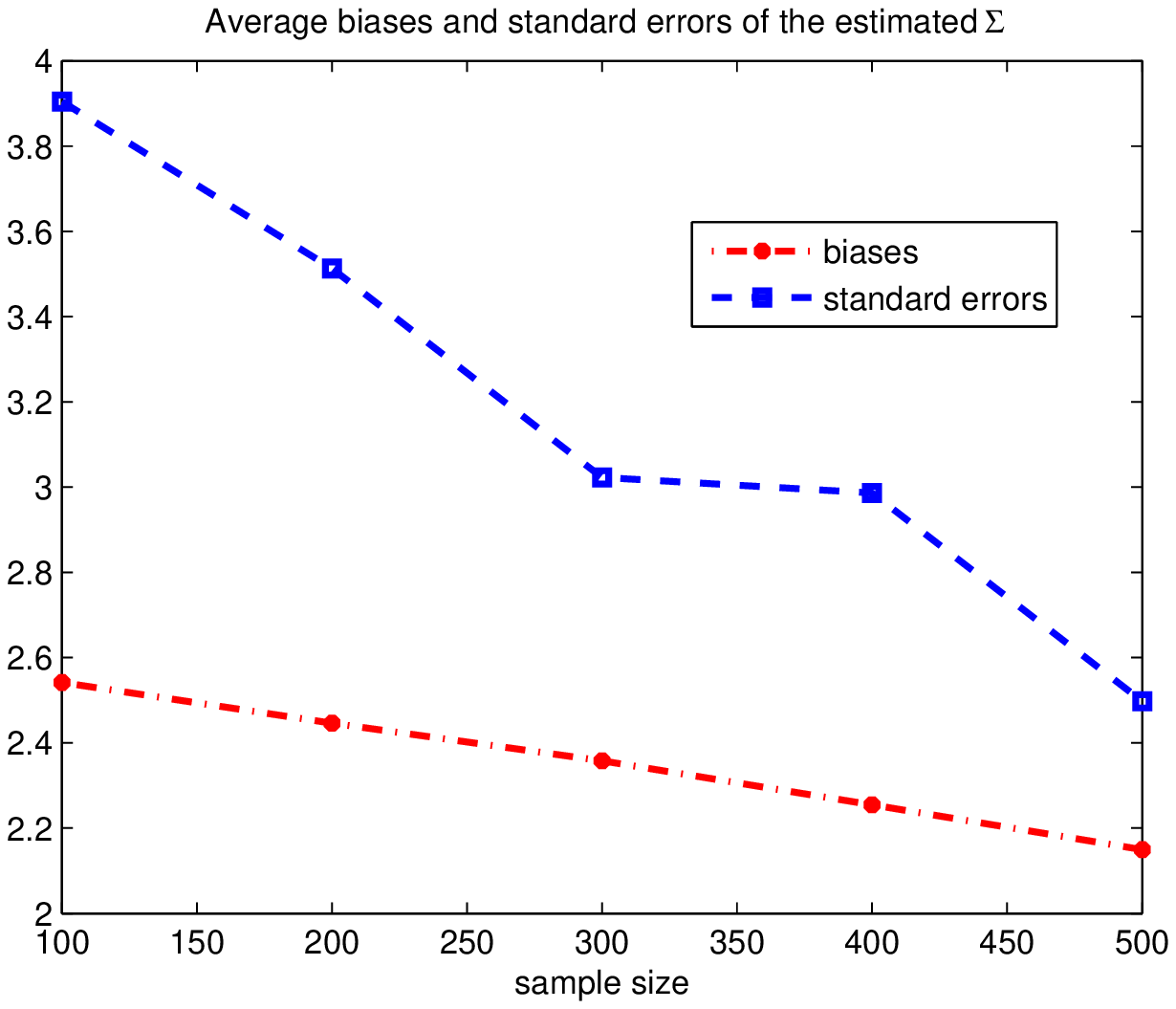}
\caption{Average bias and standard error of the MCE's for $a_0$ (top left), $\mu$ (top right), and $\Sigma$ (bottom), as a function of the sample size $n$.}
\label{fig:results_simu}
\end{figure}

Figure \ref{fig:sample_simu} shows one random sample of 100 observations generated from the model. We show the average biases and standard errors of the estimates as functions of the sample size in Figure \ref{fig:results_simu} . Here, the average bias and standard error of the estimates of $\Sigma$ are the $L_2$ norms of the average bias and standard error matrices, respectively.  As expected from Corollary \ref{coro:consist} and Theorem \ref{thm:clt}, both the bias and the standard error reduce to 0 as sample size grows to infinity. The numerical results are summarized in Table \ref{tab:mc_1}.

Finally, we point out that the choice of the region of integration
$S$ is important. A larger $S$ usually leads to more accurate
estimates, but could also result in more computational complexity.
We do not investigate this issue in this paper. However, based on
our simulation experience, an $S$ that covers most of the points
$(a,b)\in\mathbb{R}^2$ such that $[a,b]$ hits some of the observed
intervals, is a good choice as a rule of thumb. In our simulation,
$E(A)\approx [20, 40]$, by ignoring the small probability
$P(\eta<0)$. Therefore, we choose
$S=\left\{(x-y,x+y): 20\leq x\leq 40, 0\leq y\leq 10\right\}$, and the estimates are satisfactory. \\

\begin{table}[ht]
\center
\caption{Average biases and standard errors of the MCE's of the model parameters in the simulation study.}
\label{tab:mc_1}
\bigskip
\begin{tabular}{rrrrrrrrr}
\hline
           &            &            &            &     {\bf } &     {\bf } &            &            &            \\

         n & \multicolumn{ 2}{c}{$a_0$=1} &       & \multicolumn{ 2}{c}{$\mu$=20} &      & \multicolumn{ 2}{c}{$\Sigma$} \\
\hline
           &            &            &            &            &            &            &            &            \\

           & {\bf bias} &  {\bf ste} &     {\bf } & {\bf bias} &  {\bf ste} &     {\bf } & {\bf bias} &  {\bf ste} \\

           &            &            &            &            &            &            &            &            \\

       100 &     0.0683 &     0.1289 &            &     1.1648 &     1.7784 &            &     4.1166 &     5.7951 \\

       200 &     0.0387 &     0.0457 &            &     0.4569 &     0.5924 &            &     3.8581 &     4.0558 \\

       300 &     0.0274 &     0.0326 &            &     0.1831 &     0.2598 &            &     3.0317 &     3.9042 \\

       400 &     0.0157 &     0.0227 &            &     0.1575 &     0.2044 &            &     2.8210 &     3.5128 \\

       500 &     0.0128 &     0.0161 &            &     0.1197 &     0.1790 &            &     2.1494 &     2.4973 \\
\hline
\end{tabular}
\end{table}

\section{A real data application}\label{sec:real}
\begin{figure}[ht]
\centering
\includegraphics[ height=1.500in, width=2.200in]{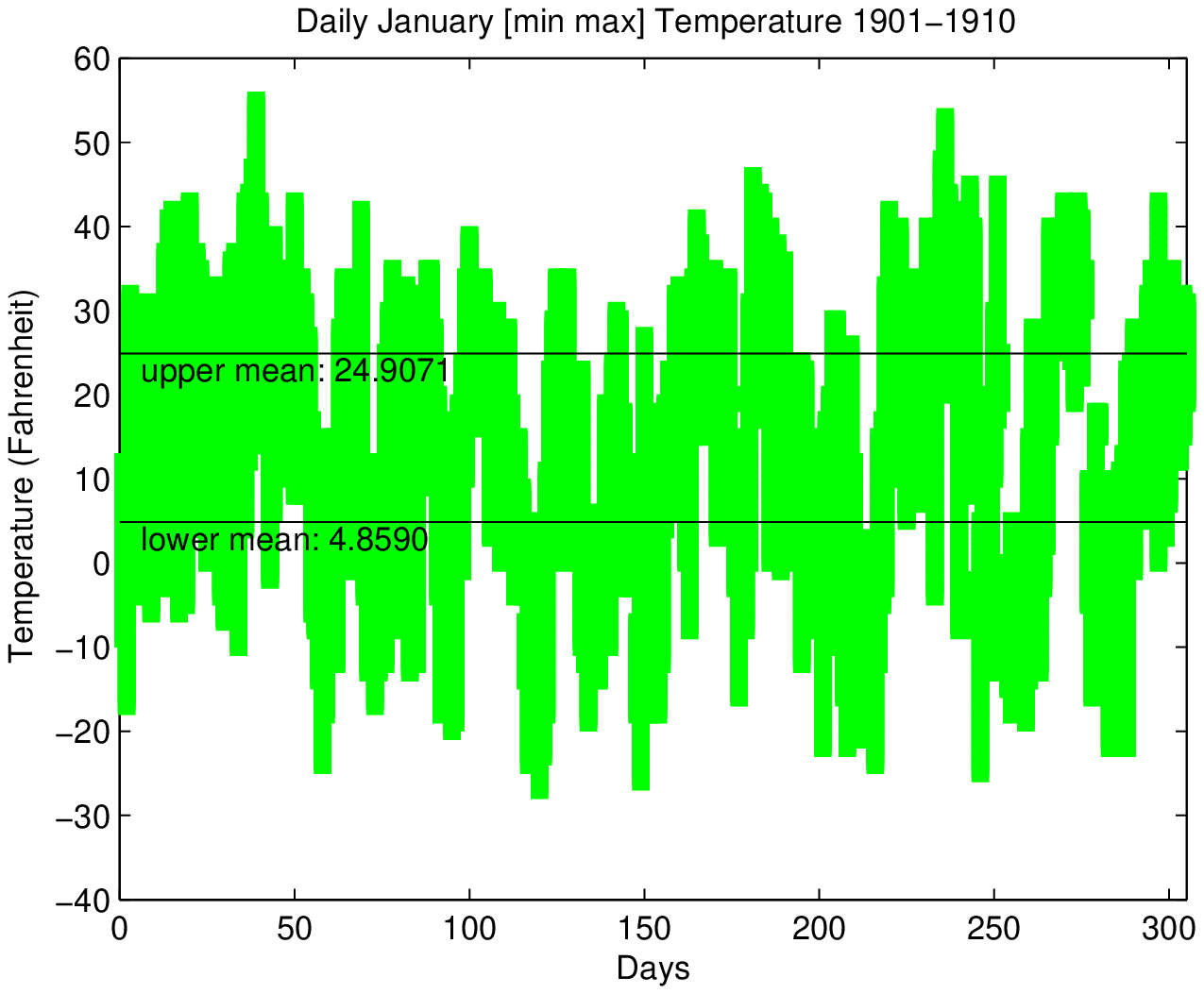}
\includegraphics[ height=1.500in, width=2.200in]{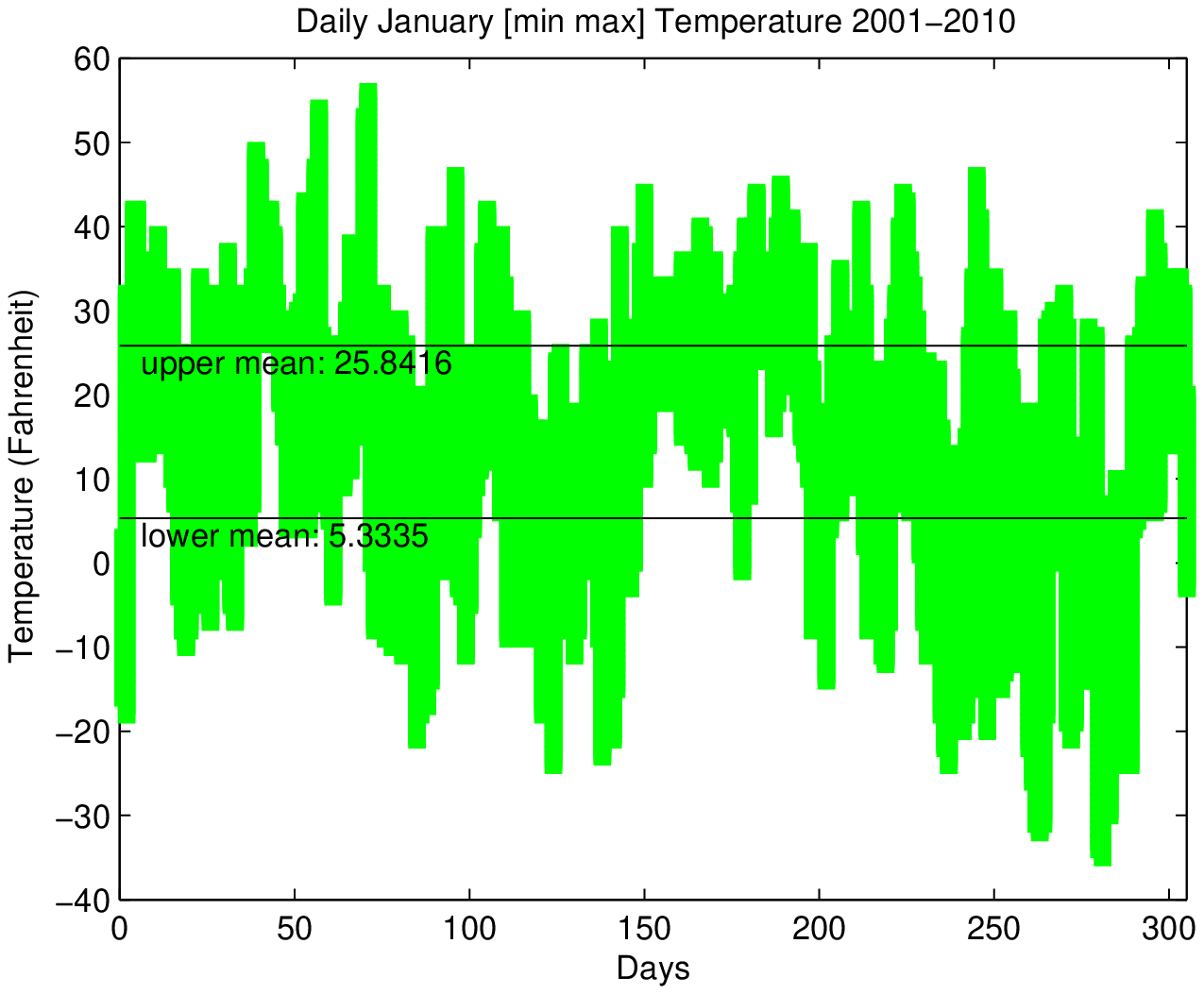}\\
\caption{Plots of daily January temperature range 1901-1910 (left) and 2001-2010 (right). On each plot, the model fitted mean is the interval between the two horizontal lines, and the moment estimate of mean is the interval between the two dashed horizontal lines.  }
\label{fig:real}
\end{figure}

In this section, we apply our Normal hierarchical model and minimum contrast
estimator to analyze the daily temperature range data. We consider two data
sets containing ten years of daily minimum and maximum temperatures in January,
in Granite Falls, Minnesota (latitude 44.81241, longitude 95.51389) from 1901
to 1910, and from 2001 to 2010, respectively. Each data set, therefore, is
constituted of 310 observations of the form: [minimum temperature, maximum
temperature] . We obtained these data from the National Weather Service, and
all observations are in Fahrenheit. The plot of the data is shown in
Figure \ref{fig:real}. The obvious correlations of the data play no roles here. \\

Same as in the simulation, we assume a bivariate normal distribution for $(\epsilon, \eta)$ and $I_0=[a_0, a_0+1]$ has length 1. The initial parameter values are computed according to (\ref{mm-1})-(\ref{mm-3}), and the weight function $W\equiv 1$. The minimum contrast estimates for the model parameters are:
\begin{itemize}
  \item Data set 1 (1901-1910):
    \begin{equation*}
    \hat{a}_{0,1}=0.2495, \hat{\mu}_1=19.8573,
    \hat{\Sigma}_1=\begin{bmatrix}207.1454 & -44.8547\\-44.8547 & 102.5263\end{bmatrix},
    \end{equation*}

  \item Data set 2 (2001-2010):
    \begin{equation*}
    \hat{a}_{0,2}=0.2614, \hat{\mu}_2=20.4722,
    \hat{\Sigma}_2=\begin{bmatrix}318.9283 & -84.0892\\-84.0892 & 68.4783\end{bmatrix}.
    \end{equation*}
 \end{itemize}

\begin{figure}[ht]
\centering
\includegraphics[ height=1.500in, width=2.200in]{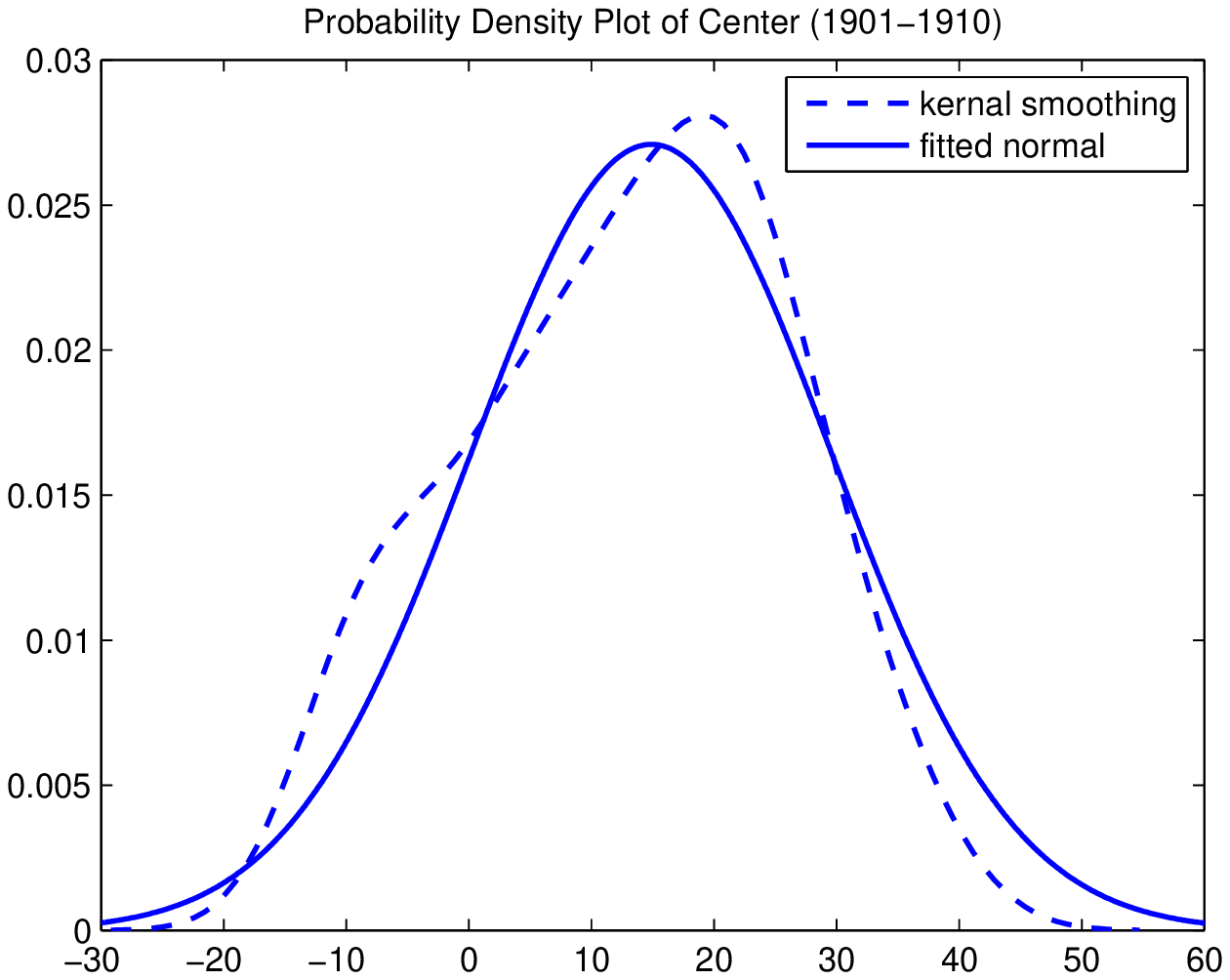}
\includegraphics[ height=1.500in, width=2.200in]{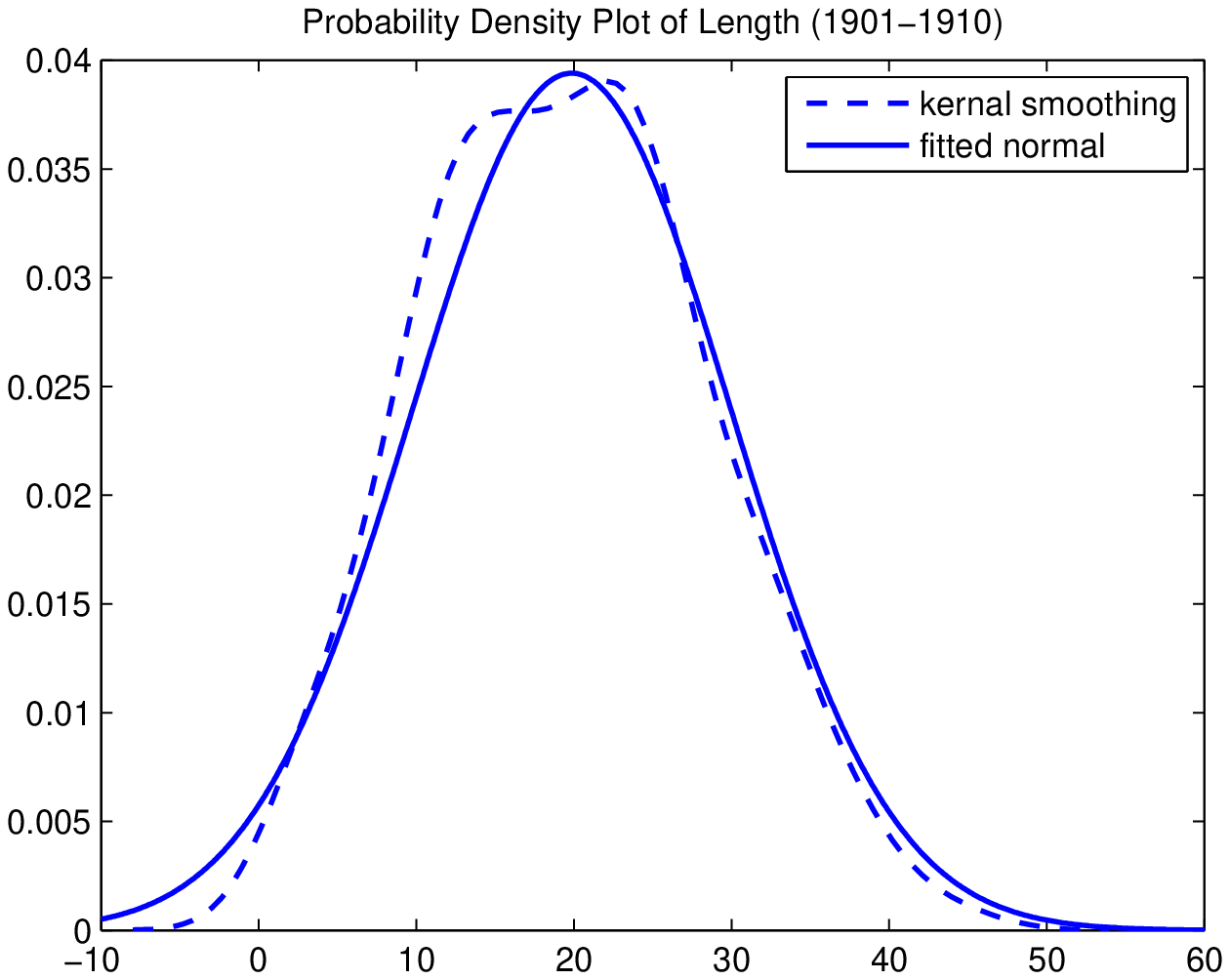}\\
\includegraphics[ height=1.500in, width=2.200in]{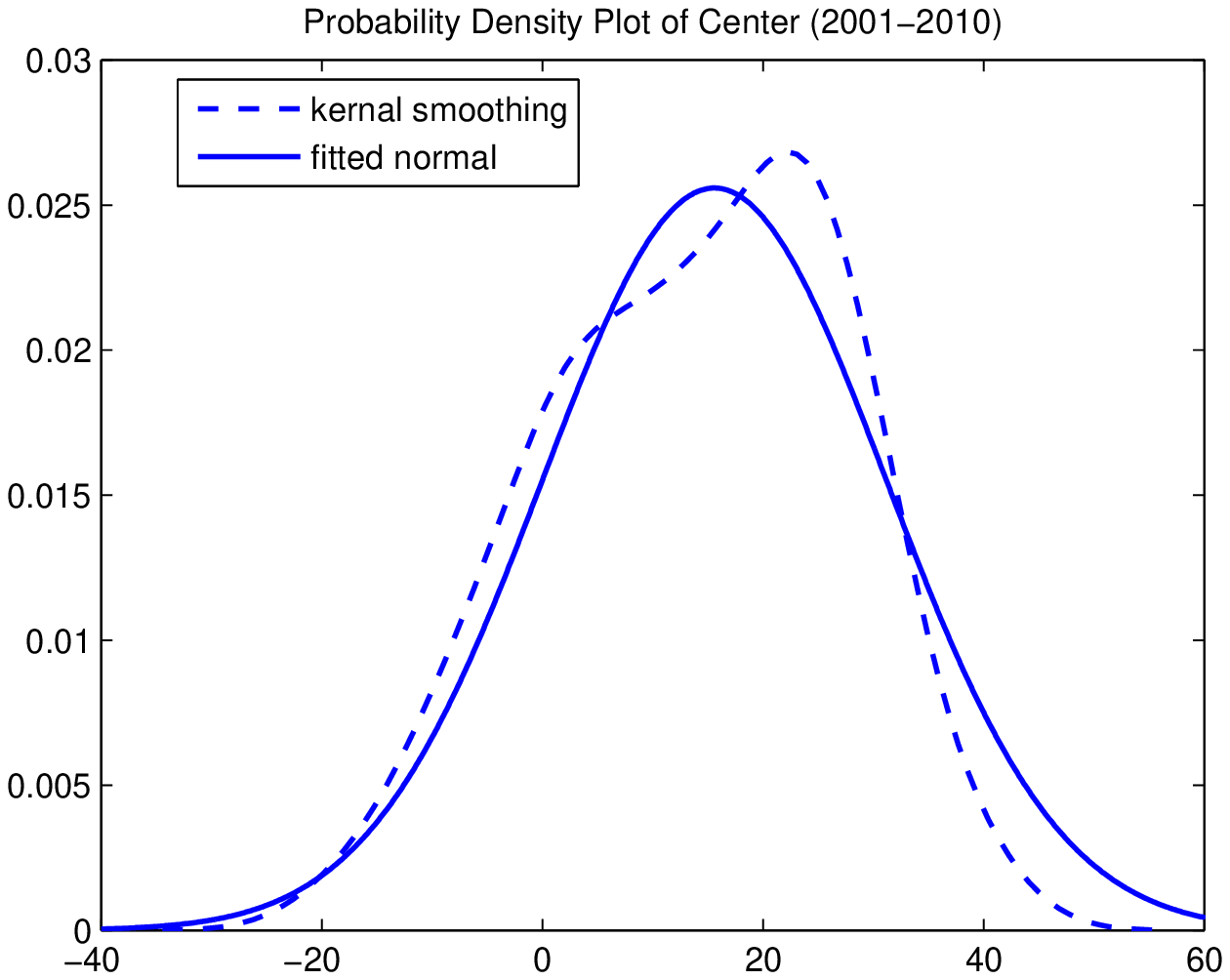}
\includegraphics[ height=1.500in, width=2.200in]{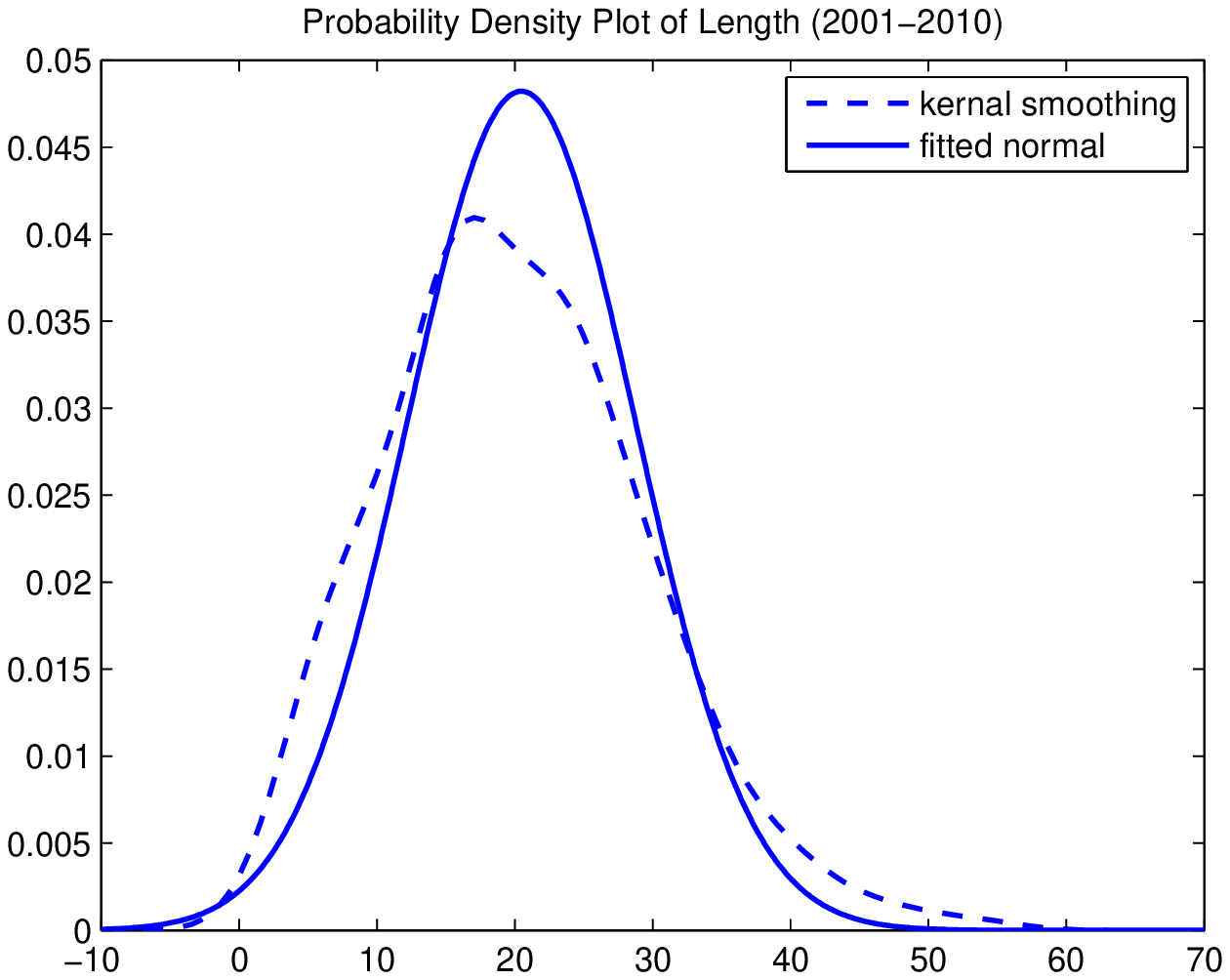}\\
\caption{Plots of the kernel smoothing density and the fitted Normal probability density for the centers and the lengths of the two data sets. }
\label{fig:pdf_plot}
\end{figure}

Recall that the center and the length of the Normal hierarchical random interval are $\epsilon+(a_0+\frac{1}{2})\eta$ and $|\eta|$($\approx\eta$ for the two considered data sets), respectively. Therefore, they are assumed to follow Normal distributions with means $(a_0+\frac{1}{2})\mu$ and $\mu$, and variances $\sigma_1^2+\left(a_0+\frac{1}{2}\right)^2\sigma_2^2+\left(2a_0+1\right)\sigma_{12}^2$ and $\sigma_2^2$, respectively. To assess the goodness-of-fit, we compare the fitted Normal distributions with the corresponding empirical distributions for both the center and the length of the two data sets. The results are shown in Figure \ref{fig:pdf_plot}. For the interval length of data 2 (2001-2010), the fitted Normal distribution is slightly more deviated from the empirical distribution, due to the skewness and heavy tail of the data. All the other three plots show very good fittings of our model to the data.\\

Denote by $A_1$ and $A_2$ respectively the random intervals from which the two data sets are drawn. The model fitted mean and variance for $A_1$ and $A_2$ are found to be:
\begin{eqnarray*}
&&\hat{\text{E}}(A_1)=\left[4.8590, 24.9071\right], \widehat{\text{Var}}(A_1)=221.2313;\\
&&\hat{\text{E}}(A_2)=\left[5.3335, 25.8416\right],
\widehat{\text{Var}}(A_2)=247.3275.
\end{eqnarray*}
Both mean and variance of the recent data are larger than those of the data 100 years ago. The two model fitted means are also shown on the data plots blue as the intervals between the solid horizontal lines in Figure \ref{fig:real}. In addition, the correlation coefficient of $(\epsilon, \eta)$ is $-0.3078$ for data set 1 and $-0.5690$ for data set 2, suggesting a negative correlation between the location and the length for the January temperature range data in general. That is, colder days tend to have larger temperature ranges, and, this relationship is stronger in the more recent data. \\

Finally, we point out that some of the parameters can
be easily estimated by simple traditional methods. For example, by
averaging the two interval ends respectively, we get the moment
estimates for the two means:
\begin{eqnarray*}
&&\hat{\text{E}}_{M}(A_1)=\left[3.5323, 22.1968\right],\\
&&\hat{\text{E}}_{M}(A_2)=\left[3.8323, 23.6903\right].
\end{eqnarray*}
They are shown in Figure \ref{fig:real} as the intervals between the dashed horizontal lines, in comparison with our model fitted means. Further, the sample correlations between the interval centers and lengths are computed as $-0.1502$ and $-0.3148$ for data sets 1 and 2, respectively. These estimates can be viewed as a preliminary analysis. Our model and the MCE of the parameters refine it and provide a more systematic understanding of the data, by examining their geometric structure in the framework of random sets.

\section{Conclusion}\label{sec:conclu}
In this paper we introduced a new model of random sets (specifically for random intervals). In many practical situations data are not completely known, or are only known with some margins of error, and it is a very important issue to consider a model which extends normality for ordinary (numerical) data. Our hierarchical normal model extends normality for point-valued random variables, and is quite flexible in the sense that it is well suited for both theoretical investigations and for simulations and real data analysis. To these goals we have defined a minimum contrast estimator for the model parameters, and we have proved its consistency and asymptotic normality. We carry out simulation experiments, and, finally we apply our model to a real data set (daily temperature range data obtained from the National Weather Service). Our approach is suitable for extensions to models in higher dimensions, e.g., a factor model for multiple random intervals, or more general random sets, including possible extensions to spherical random sets.

\section{Proofs}\label{sec:proofs}
\subsection{Proof of Theorem \ref{thm:strong-consist}}
Assume by contradiction that $\hat{\boldsymbol{\theta}}_n$ does not converge to $\boldsymbol{\theta}_0$ almost surely. Then, there exists an $\epsilon>0$ such that
\begin{equation}
  P(\left\{\omega:
  \limsup_{n\rightarrow\infty} \left\|\hat{\boldsymbol{\theta}}_n(\omega)-\boldsymbol{\theta}_0\right\|
  \geq\epsilon\right\})>0. \nonumber
\end{equation}
Let $F:=\left\{\omega:
  \limsup_{n\rightarrow\infty} \left\|\hat{\boldsymbol{\theta}}_n(\omega)-\boldsymbol{\theta}_0\right\|
  \geq\epsilon\right\}$ and $\Lambda:=\Theta\cap\left\{\boldsymbol{\theta}:\left\|\boldsymbol{\theta}-
  \boldsymbol{\theta}_0\right\|\geq\epsilon\right\}$. By the compactness of $\Lambda$, for every $\omega\in F$, there exists a convergent subsequence $\left\{\hat{\boldsymbol{\theta}}_{n_i}(\omega)\right\}$ of $\left\{\hat{\boldsymbol{\theta}}_{n}(\omega)\right\}$ such that
\begin{equation}
  \hat{\boldsymbol{\theta}}_{n_i}(\omega)\rightarrow \tilde{\boldsymbol{\theta}}\left(\omega\right) \in\Lambda, \nonumber
\end{equation}
as $i\rightarrow\infty$. Since $\bdsm{\theta}_0$ is the true underlying parameter vector that generates $X(n)$, from Definition \ref{def:cf}, $M\left(X(n);\bdsm{\theta}_0\right)$ converges to $N\left(\bdsm{\theta}_0, \bdsm{\theta}_0\right)$ almost surely, and any subsequence converges too. So we have
\begin{equation*}
  \lim_{i\rightarrow\infty}M(X(n_i);\boldsymbol{\theta}_0)
  =N(\boldsymbol{\theta}_0,\boldsymbol{\theta}_0).
\end{equation*}
On the other hand, almost surely,
\begin{eqnarray}
 &&\lim_{i\rightarrow\infty}M(X(n_i);\boldsymbol{\theta}_0)\\
  &=&\liminf_{i\rightarrow\infty}M(X(n_i);\boldsymbol{\theta}_0)\nonumber\\
  &\geq&\liminf_{i\rightarrow\infty}M(X(n_i);\hat{\boldsymbol{\theta}}_{n_i})\nonumber\\
  &=&\liminf_{i\rightarrow\infty}\left\{M(X(n_i);\hat{\boldsymbol{\theta}}_{n_i})-M(X(n_i);\tilde{\boldsymbol{\theta}})
  +M(X(n_i);\tilde{\boldsymbol{\theta}})\right\}\nonumber\\
  &\geq&\liminf_{i\rightarrow\infty}\left\{M(X(n_i);\hat{\boldsymbol{\theta}}_{n_i})
  -M(X(n_i);\tilde{\boldsymbol{\theta}})\right\}
  +\liminf_{i\rightarrow\infty}\left\{M(X(n_i);\tilde{\boldsymbol{\theta}})\right\}\nonumber\\
  &=&\liminf_{i\rightarrow\infty}\left\{M(X(n_i);\tilde{\boldsymbol{\theta}})\right\}\label{equicon}\\
  &=&\lim_{i\rightarrow\infty}\left\{M(X(n_i);\tilde{\boldsymbol{\theta}})\right\}\nonumber\\
  &=&N(\tilde{\boldsymbol{\theta}};\boldsymbol{\theta}_0).\nonumber
\end{eqnarray}
Equation (\ref{equicon}) follows from the equicontinuity of $M(X(n);\boldsymbol{\theta})$.

Therefore,
\begin{equation}\label{contra}
  P(\left\{\omega: N(\boldsymbol{\theta}_0,\boldsymbol{\theta}_0)
  \geq N(\tilde{\boldsymbol{\theta}}(\omega),\boldsymbol{\theta}_0)\right\})>0,
\end{equation}
where $\tilde{\boldsymbol{\theta}}(\omega)\in\Lambda$ and consequently $\tilde{\boldsymbol{\theta}}\neq\boldsymbol{\theta}_0$. But from the assumptions, $N(\boldsymbol{\theta}_0,\boldsymbol{\theta}_0)
  < N(\tilde{\boldsymbol{\theta}}(\omega),\boldsymbol{\theta}_0), \forall \omega$. This contradicts (\ref{contra}). Hence the desired result follows.

\subsection{Proof of Theorem \ref{thm:clt}}
From Taylor's Theorem, we have
\begin{eqnarray}
 0 &=&\frac{\partial H}{\partial\theta_i}\left(X\left(n\right);\bdsm{\theta}_n^H\right)\nonumber\\
 &=&\frac{\partial H}{\partial\theta_i}\left(X\left(n\right);\bdsm{\theta}_0\right)
 +\sum\limits_{j=1}^{p}\left(\theta^H_{n,j}-\theta_{0,j}\right)
 \frac{\partial^2 H}{\partial\theta_j\partial\theta_i}\left(X\left(n\right);\bdsm{\theta}_0\right)\nonumber\\
 &&+\frac{1}{2}\left[\sum\limits_{j=1}^{p}\left(\theta^H_{n,j}-\theta_{0,j}\right)\frac{\partial}{\partial\theta_j}\right]^2
 \frac{\partial H}{\partial\theta_i}\left(X\left(n\right);\bdsm{\epsilon}_n\right)\nonumber\\
 &=&\frac{\partial H}{\partial\theta_i}\left(X\left(n\right);\bdsm{\theta}_0\right)\nonumber\\
 &&+\sum\limits_{j=1}^{p}\left(\theta^H_{n,j}-\theta_{0,j}\right)\left[
 \frac{\partial^2H}{\partial\theta_j\partial\theta_i}\left(X(n);\bdsm{\theta}_0\right)+\frac{1}{2}
 \sum\limits_{l=1}^{p}\left(\theta^H_{n,l}-\theta_{0,l}\right)\frac{\partial^3H}
 {\partial\theta_l\partial\theta_j\partial\theta_i}\left(X\left(n\right);\bdsm{\epsilon}_n\right)
 \right],\nonumber
\end{eqnarray}
for $i=1,\cdots,p$, where $\bdsm{\epsilon}_n$ lies between $\bdsm{\theta}_0$ and $\bdsm{\theta}_n^H$. Writing the above equations in matrix form, we get
\begin{eqnarray}
 &&\frac{\partial H}{\partial\bdsm{\theta}}\left(X(n);\bdsm{\theta}_0\right)\nonumber\\
 &&+\left[\frac{\partial^2H}{\partial\bdsm{\theta}^2}
 \left(X(n);\bdsm{\theta}_0\right)+\frac{1}{2}\sum\limits_{j=1}^{p}\left(\theta^H_{n,j}-\theta_{0,j}\right)
 \left(\frac{\partial}{\partial\theta_j}\left(\frac{\partial^2H}{\partial\bdsm{\theta}}\right)(X(n);\bdsm{\epsilon}_n)\right)
 \right]
 \left(\bdsm{\theta}_n^H-\bdsm{\theta}_0\right)\nonumber\\
 &&=0\label{thm3:eqn1}.
\end{eqnarray}
Observe, by taking derivatives under the integral sign, that $\forall i,j$,
\begin{eqnarray}
 &&\frac{\partial^2H}{\partial\theta_j\partial\theta_i}\left(X(n);\bdsm{\theta}_0\right)\nonumber\\
 &=&\frac{\partial^2H}{\partial\theta_j\partial\theta_i}\iint\limits_{S}
 \left[T_{\boldsymbol{\theta}}([a,b])-\hat{T}([a,b];X(n))\right]^2W(a,b)\mathrm{d}a\mathrm{d}b,\nonumber\\
 &=&\frac{\partial}{\partial\theta_j}2\iint\limits_{S}\left[T_{\boldsymbol{\theta}}([a,b])-\hat{T}([a,b];X(n))\right]
 \frac{\partial T_{\bdsm{\theta}_0}}{\partial\theta_i}([a,b])W(a,b)\mathrm{d}a\mathrm{d}b,\nonumber\\
 &=&2\iint\limits_{S}\left[T_{\boldsymbol{\theta}}([a,b])-\hat{T}([a,b];X(n))\right]
 \frac{\partial^2T_{\bdsm{\theta}_0}}{\partial\theta_j\partial\theta_i}([a,b])W(a,b)\mathrm{d}a\mathrm{d}b\nonumber\\
 &&+2\iint\limits_{S}\left(\frac{\partial T_{\bdsm{\theta}_0}}{\partial\theta_j}
 \frac{\partial T_{\bdsm{\theta}_0}}{\partial\theta_i}\right)([a,b])W(a,b)\mathrm{d}a\mathrm{d}b\nonumber\\
 &:=&I+II.\nonumber
\end{eqnarray}
The first term is
\begin{eqnarray}
 I&=&2\iint\limits_{S}\left(T_{\bdsm{\theta}_0}\left([a,b]\right)-\frac{1}{n}\sum_{k=1}^{n}Y_k\left(a,b\right)\right)
  \frac{\partial^2 T_{\bdsm{\theta}_0}}{\partial\theta_j\partial\theta_i}([a,b])W(a,b)\mathrm{d}a\mathrm{d}b\nonumber\\
  &=&\frac{2}{n}\sum_{k=1}^{n}\iint\limits_{S}\left[T_{\bdsm{\theta}_0}\left([a,b]\right)-Y_k\left(a,b\right)\right]
  \frac{\partial^2 T_{\bdsm{\theta}_0}}{\partial\theta_j\partial\theta_i}([a,b])W(a,b)\mathrm{d}a\mathrm{d}b\nonumber\\
  &=&o_P(1),\nonumber
\end{eqnarray}
according to the strong law of large numbers for i.i.d. random variables. Therefore,
\begin{equation}
 \frac{\partial^2H}{\partial\theta_j\partial\theta_i}\left(X(n);\bdsm{\theta}_0\right)
 =o_P(1)+2\iint\limits_{S}\left(\frac{\partial T_{\bdsm{\theta}_0}}{\partial\theta_j}
 \frac{\partial T_{\bdsm{\theta}_0}}{\partial\theta_i}\right)([a,b])W(a,b)\mathrm{d}a\mathrm{d}b,\nonumber
\end{equation}
$\forall i,j$. In matrix form,
\begin{equation}\label{thm3:eqn2}
 \frac{\partial^2H}{\partial\bdsm{\theta}^2}\left(X(n);\bdsm{\theta}_0\right)
 =o_P(1)+2\iint\limits_{S}\left(\frac{\partial T_{\bdsm{\theta}_0}}{\partial\bdsm{\theta}}\right)
 \left(\frac{\partial T_{\bdsm{\theta}_0}}{\partial\bdsm{\theta}}\right)^{T}([a,b])W(a,b)\mathrm{d}a\mathrm{d}b.
\end{equation}
Observe again that $\forall j,k,l$,
\begin{eqnarray}
 &&\left|\frac{\partial^3H(X(n);\bdsm{\epsilon}_n)}{\partial\theta_j\partial\theta_k\partial\theta_l}\right|\nonumber\\
 &\leq&2\iint\limits_{S}\left|\left[T_{\bdsm{\epsilon}_n}([a,b])-\hat{T}([a,b];X(n))\right]
 \frac{\partial^3T_{\bdsm{\epsilon}_n}}{\partial\theta_j\partial\theta_k\partial\theta_l}
 ([a,b])W(a,b)\mathrm{d}a\mathrm{d}b\right|\nonumber\\
 &&+2\left|\iint\limits_{S}\left[\left(\frac{\partial T_{\bdsm{\epsilon}_n}}{\partial\theta_j}
 \frac{\partial^2T_{\bdsm{\epsilon}_n}}{\partial\theta_k\partial\theta_l}\right)
 +\left(\frac{\partial^2T_{\bdsm{\epsilon}_n}}{\partial\theta_j\partial\theta_k}\frac{\partial T_{\bdsm{\epsilon}_n}}
 {\partial\theta_l}\right)
 +\left(\frac{\partial^2T_{\bdsm{\epsilon}_n}}{\partial\theta_j\partial\theta_l}\frac{\partial T_{\bdsm{\epsilon}_n}}
 {\partial\theta_k}\right)\right]([a,b])W(a,b)\mathrm{d}a\mathrm{d}b\right|\nonumber\\
 &\leq&4\iint\limits_{S}\left|\frac{\partial^3T_{\bdsm{\epsilon}_n}}{\partial\theta_j\partial\theta_k\partial\theta_l}
 ([a,b])W(a,b)\mathrm{d}a\mathrm{d}b\right|\nonumber\\
 &&+2\left|\iint\limits_{S}\left[\left(\frac{\partial T_{\bdsm{\epsilon}_n}}{\partial\theta_j}
 \frac{\partial^2T_{\bdsm{\epsilon}_n}}{\partial\theta_k\partial\theta_l}\right)
 +\left(\frac{\partial^2T_{\bdsm{\epsilon}_n}}{\partial\theta_j\partial\theta_k}\frac{\partial T_{\bdsm{\epsilon}_n}}
 {\partial\theta_l}\right)
 +\left(\frac{\partial^2T_{\bdsm{\epsilon}_n}}{\partial\theta_j\partial\theta_l}\frac{\partial T_{\bdsm{\epsilon}_n}}
 {\partial\theta_k}\right)\right]([a,b])W(a,b)\mathrm{d}a\mathrm{d}b\right|\nonumber\\
 &:=&C_1(\bdsm{\epsilon}_n)\leq C_2,\nonumber
\end{eqnarray}
$\forall\bdsm{\epsilon}_n\in\Theta$, by the compactness of $\Theta$. This, together with the strong consistency of $\bdsm{\theta}_n^H$, gives
\begin{eqnarray}
 &&\frac{1}{2}\sum\limits_{j=1}^{p}\left(\theta^H_{n,j}-\theta_{0,j}\right)
 \left(\frac{\partial}{\partial\theta_j}\left(\frac{\partial^2H}{\partial\theta_k\partial\theta_l}\right)
 (X(n);\bdsm{\epsilon}_n)\right)\nonumber\\
 &=&\frac{1}{2}\sum\limits_{j=1}^{p}o_P(1)\frac{\partial^3H(X(n);\bdsm{\epsilon}_n)}
 {\partial\theta_j\partial\theta_k\partial\theta_l}\nonumber\\
 &=&o_P(1),\nonumber
\end{eqnarray}
$\forall k,l$. Equivalently, in matrix form,
\begin{equation}\label{thm3:eqn3}
 \frac{1}{2}\sum\limits_{j=1}^{p}\left(\theta^H_{n,j}-\theta_{0,j}\right)
 \left(\frac{\partial}{\partial\theta_j}\left(\frac{\partial^2H}{\partial\bdsm{\theta}}\right)
 (X(n);\bdsm{\epsilon}_n)\right)=o_P(1).
\end{equation}
By the multivariate Slutsky's theorem, Proposition \ref{prop:parH}, together with equation (\ref{thm3:eqn1}), (\ref{thm3:eqn2}), and (\ref{thm3:eqn3}), yields the desired result.\\



\newpage
\section{Appendix}
\subsection{Proof of Proposition \ref{prop:cf}}
Notice that $\hat{T}([a,b];X(n))$ is the sample mean of i.i.d. random variables $Y_i:\Omega\rightarrow \mathbb{R}$ defined as:
\begin{equation}\label{Y_def}
  Y_i=\begin{cases}
    1, & \text{if} \ X_i\cap[a,b]\neq\emptyset, \\
    0, & \text{otherwise}.
  \end{cases}.
\end{equation}
Therefore, an application of the strong law of large numbers in the classical case yields:
\begin{equation*}
  \frac{1}{n}\sum_{i=1}^{n}Y_i
  \stackrel{a.s.}{\rightarrow} EY_1=P\left(X_1\cap[a,b]\neq\emptyset\right)
  =T_{\boldsymbol{\theta}_0}\left([a,b]\right),\ \text{as}\ n\to\infty,
\end{equation*}
$\forall a,b: -\infty<a\leq b<\infty$, and assuming $\boldsymbol{\theta}_0$ is the true parameter value. That is,
\begin{equation}
  \hat{T}\left([a,b];X(n)\right)\stackrel{a.s.}{\rightarrow}T_{\boldsymbol{\theta}_0}\left([a,b]\right), \nonumber
\end{equation}
as $n\to\infty$. It follows immediately that
\begin{equation}
  \left[\hat{T}([a,b];X(n))-T_{\boldsymbol{\theta}_0}\left([a,b]\right)\right]^2W(a,b)\stackrel{a.s.}{\rightarrow}0. \nonumber
\end{equation}
Notice that $\forall a,b: -\infty<a\leq b<\infty$, $\left[\hat{T}([a,b];X(n))-T_{\boldsymbol{\theta}_0}\left([a,b]\right)\right]^2W(a,b)$ is uniformly bounded by $4C$. By the bounded convergence theorem,
\begin{equation}
  \iint\limits_{S}\left[\hat{T}([a,b];X(n))-T_{\boldsymbol{\theta}_0}\left([a,b]\right)\right]^2W(a,b)\mathrm{d}a\mathrm{d}b
  \stackrel{a.s.}{\rightarrow}\iint\limits_{S}0\cdot \mathrm{d}a\mathrm{d}b=0, \nonumber
\end{equation}
given any $S\subset\mathbb{R}^2$ with finite Lebesgue measure. This verifies that
\begin{equation}\label{eqn:N1}
  P_{\boldsymbol{\theta}}\left\{\omega: \lim_{n\to\infty}H\left(X(n);\boldsymbol{\theta}\right)=0\right\}=1.
\end{equation}
Similarly, we also get
\begin{equation}\label{eqn:N2}
  P_{\boldsymbol{\theta}}\left\{\omega: \lim_{n\to\infty}H\left(X(n);\boldsymbol{\zeta}\right)=
  \iint\limits_{S}\left[T_{\boldsymbol{\theta}}([a,b])-T_{\boldsymbol{\zeta}}([a,b])\right]^2W(a,b)\mathrm{d}a\mathrm{d}b\right\}=1,
\end{equation}
$\forall \boldsymbol{\theta},\boldsymbol{\zeta}\in\Theta$. Equations (\ref{eqn:N1}) and (\ref{eqn:N2}) together imply
\begin{equation}\label{eqn:N}
  N(\boldsymbol{\theta},\boldsymbol{\zeta})=
  \iint\limits_{S}\left[T_{\boldsymbol{\theta}}([a,b])-T_{\boldsymbol{\zeta}}([a,b])\right]^2W(a,b)\mathrm{d}a\mathrm{d}b,\
  \boldsymbol{\theta}, \boldsymbol{\zeta}\in\Theta.
\end{equation}
By Assumption \ref{aspt:2}, $T_{\boldsymbol{\theta}}([a,b])\neq T_{\boldsymbol{\zeta}}([a,b])$, for $\boldsymbol{\theta}\neq\boldsymbol{\zeta}$, except on a Lebesgue set of measure 0. This together with (\ref{eqn:N}) gives
\begin{equation}
  N(\boldsymbol{\theta},\boldsymbol{\theta})<N(\boldsymbol{\theta},\boldsymbol{\zeta}),\
  \forall\ \boldsymbol{\theta}\neq\boldsymbol{\zeta},\ \boldsymbol{\theta},\boldsymbol{\zeta}\in\Theta,\nonumber
\end{equation}
which proves that $H(X(n);\boldsymbol{\theta})$, $\theta\in\Theta$ is a family of contrast functions. To see the equicontinuity of $H(X(n);\boldsymbol{\theta})$, notice that $\forall\boldsymbol{\theta}_1,\boldsymbol{\theta}_2\in\Theta$, we have
\begin{eqnarray*}
 &&\left|H(X(n);\boldsymbol{\theta}_1)-H(X(n);\boldsymbol{\theta}_2)\right|\\
  &=&|\iint\limits_{S}\left(T_{\boldsymbol{\theta}_1}([a,b])-\hat{T}([a,b];X(n))\right)^2W(a,b)\mathrm{d}a\mathrm{d}b\\
  &&-\iint\limits_{S}\left(T_{\boldsymbol{\theta}_2}([a,b])-\hat{T}([a,b];X(n))\right)^2W(a,b)\mathrm{d}a\mathrm{d}b|\\
  &=&|\iint\limits_{S}\left(T_{\boldsymbol{\theta}_1}([a,b])-T_{\boldsymbol{\theta}_2}([a,b])\right)
  \left(T_{\boldsymbol{\theta}_1}([a,b])+T_{\boldsymbol{\theta}_2}([a,b])
  -2\hat{T}([a,b];X(n))\right)W(a,b)\mathrm{d}a\mathrm{d}b|\\
  &\leq&4C\iint\limits_{S}\left|T_{\boldsymbol{\theta}_1}([a,b])-T_{\boldsymbol{\theta}_2}([a,b])\right|\mathrm{d}a\mathrm{d}b,
\end{eqnarray*}
since, by definition (\ref{H_def}), $|W(a,b)|$ is uniformly bounded by $C$, $\forall a,b: -\infty<a\leq b.$
Then the equicontinuity of $H(X(n);\boldsymbol{\theta})$ follows from the continuity of $T_{\boldsymbol{\theta}}([a,b])$.

\subsection{Lemma 1}
Let $H(X(n);\boldsymbol{\theta})$ be the contrast function defined in (\ref{H_def}). Under the hypothesis of Assumption \ref{aspt:4},
\begin{equation*}
  \sqrt{n}\left[\frac{\partial H}{\partial\theta_i}\left(X\left(n\right);\bdsm{\theta}_0\right)\right]
  \stackrel{\mathcal{D}}{\rightarrow}
  N\left(0,\Delta_i\right),\ \text{as}\ n\to\infty,
\end{equation*}
for $i=1,\cdots,p$, where
\begin{eqnarray}
  \Delta_i&=&4\iiiint\limits_{S\times S}\left\{P\left(X_1\cap[a,b]\neq\emptyset,X_1\cap[c,d]\neq\emptyset\right)
  -T_{\bdsm{\theta}_0}\left([a,b]\right)T_{\bdsm{\theta}_0}\left([c,d]\right)\right\}\nonumber\\
  &&\times\frac{\partial T_{\bdsm{\theta}_0}}{\partial\theta_i}\left([a,b]\right)
  \frac{\partial T_{\bdsm{\theta}_0}}{\partial\theta_i}\left([c,d]\right)
  W(a,b)W(c,d)\mathrm{d}a\mathrm{d}b\mathrm{d}c\mathrm{d}d.\nonumber
\end{eqnarray}
\begin{proof}
We will write $\frac{\partial T_{\bdsm{\theta}_0}\left([a,b]\right)}{\partial\theta_i}=
T_{\bdsm{\theta}_0}^i\left(a,b\right)$ to simplify notations. Exchanging differentiation and integration by the bounded convergence theorem, we get
\begin{eqnarray}
  &&\frac{\partial H}{\partial\theta_i}\left(X\left(n\right);\bdsm{\theta}_0\right)\label{par_H}\\
  &=&\frac{\partial}{\partial\theta_i}\iint\limits_{S}
  \left(T_{\bdsm{\theta}_0}\left([a,b]\right)-\hat{T}\left([a,b];X(n)\right)\right)^2W(a,b)\mathrm{d}a\mathrm{d}b\nonumber\\
  &=&\iint\limits_{S}\frac{\partial}{\partial\theta_i}
  \left(T_{\bdsm{\theta}_0}\left([a,b]\right)-\hat{T}\left([a,b];X(n)\right)\right)^2W(a,b)\mathrm{d}a\mathrm{d}b\nonumber\\
  &=&\iint\limits_{S}2\left(T_{\bdsm{\theta}_0}\left([a,b]\right)-\hat{T}\left([a,b];X(n)\right)\right)
  T_{\bdsm{\theta}_0}^i\left(a,b\right)W(a,b)\mathrm{d}a\mathrm{d}b.\nonumber
\end{eqnarray}
Define $Y_i(a,b)$ as in (\ref{Y_def}). Then,
\begin{eqnarray}
  (\ref{par_H})
  &=&\iint\limits_{S}2\left(T_{\bdsm{\theta}_0}\left([a,b]\right)-\frac{1}{n}\sum_{k=1}^{n}Y_k\left(a,b\right)\right)
  T_{\bdsm{\theta}_0}^i\left(a,b\right)W(a,b)\mathrm{d}a\mathrm{d}b\nonumber\\
  &=&\frac{2}{n}\iint\limits_{S}\sum_{k=1}^{n}\left(T_{\bdsm{\theta}_0}\left([a,b]\right)-Y_k\left(a,b\right)\right)
  T_{\bdsm{\theta}_0}^i\left(a,b\right)W(a,b)\mathrm{d}a\mathrm{d}b\nonumber\\
  &=&\frac{1}{n}\sum_{k=1}^{n}2\iint\limits_{S}\left(T_{\bdsm{\theta}_0}\left([a,b]\right)-Y_k\left(a,b\right)\right)
  T_{\bdsm{\theta}_0}^i\left(a,b\right)W(a,b)\mathrm{d}a\mathrm{d}b\label{eqn:parH}\\
  &:=&\frac{1}{n}\sum_{k=1}^{n}R_k.\nonumber
\end{eqnarray}
Notice that $R_k$'s are i.i.d. random variables: $\Omega\rightarrow\mathbb{R}$.\\

Let $\left\{\Delta s_1,\Delta s_2, \cdots, \Delta s_m\right\}$ be a partition of $S$, and $(a_j,b_j)$ be any point in $\Delta s_j$, $j=1,\cdots,m$. Let $\lambda=\max_{1\leq j\leq m}\left\{\text{diam}\Delta s_j\right\}$. Denote by $\Delta\sigma_j$ the area of $\Delta s_j$. By the definition of the double integral,
\begin{eqnarray}
  R_k
  &=&2\iint\limits_{S}\left(T_{\bdsm{\theta}_0}\left([a,b]\right)-Y_k\left(a,b\right)\right)
  T_{\bdsm{\theta}_0}^i\left(a,b\right)W(a,b)\mathrm{d}a\mathrm{d}b\nonumber\\
  &=&\lim_{\lambda\rightarrow 0}\left\{\sum_{j=1}^{m}\left(T_{\bdsm{\theta}_0}
  \left([a_j,b_j]\right)-Y_k\left(a_j,b_j\right)\right)
  T_{\bdsm{\theta}_0}^i\left(a_j,b_j\right)W(a_j,b_j)\Delta\sigma_j\right\}.\nonumber
\end{eqnarray}
Therefore, by the Lebesgue dominated convergence theorem,
\begin{eqnarray*}
  &&ER_k\\
  &=&2E\lim_{\lambda\rightarrow 0}\left\{\sum_{j=1}^{m}\left(T_{\bdsm{\theta}_0}
  \left([a_j,b_j]\right)-Y_k\left(a_j,b_j\right)\right)
  T_{\bdsm{\theta}_0}^i\left(a_j,b_j\right)W(a_j,b_j)\Delta\sigma_j\right\}\\
  &=&2\lim_{\lambda\rightarrow 0}\left\{\sum_{j=1}^{m}\left[E\left(T_{\bdsm{\theta}_0}
  \left([a_j,b_j]\right)-Y_k\left(a_j,b_j\right)\right)\right]
  T_{\bdsm{\theta}_0}^i\left(a_j,b_j\right)W(a_j,b_j)\Delta\sigma_j\right\}\label{eqn_1}\\
  &=&2\lim_{\lambda\rightarrow 0}\left\{\sum_{j=1}^{m}0\right\}=0.
\end{eqnarray*}

Moreover,
\begin{eqnarray*}
  &&Var(R_k)=ER_k^2\\
  &=&4E\left\{\lim_{\lambda\rightarrow 0}\left\{\sum_{j=1}^{m}\left(T_{\bdsm{\theta}_0}
  \left([a_j,b_j]\right)-Y_k\left(a_j,b_j\right)\right)
  T_{\bdsm{\theta}_0}^i\left(a_j,b_j\right)W(a_j,b_j)\Delta\sigma_j\right\}\right\}^2\\
  &=&4E\lim_{\lambda_1\rightarrow 0}\lim_{\lambda_2\rightarrow 0}
  \left\{\sum_{j_1=1}^{m_1}\left(T_{\bdsm{\theta}_0}
  \left([a_{j_1},b_{j_1}]\right)-Y_k\left(a_{j_1},b_{j_1}\right)\right)
  T_{\bdsm{\theta}_0}^i\left(a_{j_1},b_{j_1}\right)W(a_{j_1},b_{j_1})\Delta\sigma_{j_1}\right\}\\
  &&\left\{\sum_{j_2=1}^{m_2}\left(T_{\bdsm{\theta}_0}
  \left([a_{j_2},b_{j_2}]\right)-Y_k\left(a_{j_2},b_{j_2}\right)\right)
  T_{\bdsm{\theta}_0}^i\left(a_{j_2},b_{j_2}\right)W(a_{j_2},b_{j_2})\Delta\sigma_{j_2}\right\}\\
  &=&4E\lim_{\lambda_1\rightarrow 0}\lim_{\lambda_2\rightarrow 0}\sum_{j_1=1}^{m_1}\sum_{j_2=1}^{m_2}
  \left(T_{\bdsm{\theta}_0}\left([a_{j_1},b_{j_1}]\right)-Y_k\left(a_{j_1},b_{j_1}\right)\right)
  \left(T_{\bdsm{\theta}_0}\left([a_{j_2},b_{j_2}]\right)-Y_k\left(a_{j_2},b_{j_2}\right)\right)\\
  &&T_{\bdsm{\theta}_0}^i\left(a_{j_1},b_{j_1}\right)T_{\bdsm{\theta}_0}^i\left(a_{j_2},b_{j_2}\right)
  W(a_{j_1},b_{j_1})W(a_{j_2},b_{j_2})\Delta\sigma_{j_1}\Delta\sigma_{j_2}\\
  &=&4\lim_{\lambda_1\rightarrow 0}\lim_{\lambda_2\rightarrow 0}\sum_{j_1=1}^{m_1}\sum_{j_2=1}^{m_2}
  E\left(T_{\bdsm{\theta}_0}\left([a_{j_1},b_{j_1}]\right)-Y_k\left(a_{j_1},b_{j_1}\right)\right)
  \left(T_{\bdsm{\theta}_0}\left([a_{j_2},b_{j_2}]\right)-Y_k\left(a_{j_2},b_{j_2}\right)\right)\\
  &&T_{\bdsm{\theta}_0}^i\left(a_{j_1},b_{j_1}\right)T_{\bdsm{\theta}_0}^i\left(a_{j_2},b_{j_2}\right)
  W(a_{j_1},b_{j_1})W(a_{j_2},b_{j_2})\Delta\sigma_{j_1}\Delta\sigma_{j_2}\label{eqn_2}\\
  &=&4\lim_{\lambda_1\rightarrow 0}\lim_{\lambda_2\rightarrow 0}\sum_{j_1=1}^{m_1}\sum_{j_2=1}^{m_2}
  Cov\left(Y_k\left(a_{j_1},b_{j_1}\right),Y_k\left(a_{j_2},b_{j_2}\right)\right)\\
  &&T_{\bdsm{\theta}_0}^i\left(a_{j_1},b_{j_1}\right)T_{\bdsm{\theta}_0}^i\left(a_{j_2},b_{j_2}\right)
  W(a_{j_1},b_{j_1})W(a_{j_2},b_{j_2})\Delta\sigma_{j_1}\Delta\sigma_{j_2}\\
  &=&4\iiiint\limits_{S\times S}Cov\left(Y_k\left(a,b\right),Y_k\left(c,d\right)\right)
  T_{\bdsm{\theta}_0}^i\left(a,b\right)T_{\bdsm{\theta}_0}^i\left(c,d\right)
  W(a,b)W(c,d)\mathrm{d}a\mathrm{d}b\mathrm{d}c\mathrm{d}d\\
  &=&4\iiiint\limits_{S\times S}\left\{P\left(X_k\cap[a,b]\neq\emptyset,X_k\cap[c,d]\neq\emptyset\right)
  -T_{\bdsm{\theta}_0}\left([a,b]\right)T_{\bdsm{\theta}_0}\left([c,d]\right)\right\}\\
  &&T_{\bdsm{\theta}_0}^i\left(a,b\right)T_{\bdsm{\theta}_0}^i\left(c,d\right)
  W(a,b)W(c,d)\mathrm{d}a\mathrm{d}b\mathrm{d}c\mathrm{d}d.
\end{eqnarray*}
From the central limit theorem for i.i.d. random variables, the desired result follows.
\end{proof}

\subsection{Proof of Proposition \ref{prop:parH}}
By the Cram\'{e}r-Wold device, it suffices to prove
\begin{equation}\label{prop1:target}
 \sqrt{n}\sum\limits_{i=1}^{p}\lambda_i\frac{\partial H}{\partial\theta_i}\left(X(n);\bdsm{\theta}_0\right)
 \stackrel{\mathcal{D}}{\rightarrow}N\left(0,\sum\limits_{1\leq i,j\leq p}\lambda_i\lambda_j\Xi(i,j)\right),
\end{equation}
for arbitrary real numbers $\lambda_i,i=1,\cdots,p$. It is easily seen from (\ref{eqn:parH}) in the proof of Lemma 1 that
\begin{eqnarray}
 &&\sum\limits_{i=1}^{p}\lambda_i\frac{\partial H}{\partial\theta_i}\left(X(n);\bdsm{\theta}_0\right)\nonumber\\
 &=&\frac{1}{n}\sum\limits_{k=1}^{n}\left(2\sum\limits_{i=1}^{p}\lambda_i
 \iint\limits_{S}\left(T_{\bdsm{\theta}_0}\left([a,b]\right)-Y_k\left(a,b\right)\right)
 \frac{\partial T_{\bdsm{\theta}_0}}{\partial\theta_i}\left([a,b]\right)W(a,b)\mathrm{d}a\mathrm{d}b\right)\nonumber\\
 &:=&\frac{1}{n}\sum\limits_{k=1}^{n}\left(2\sum\limits_{i=1}^p\lambda_iQ_k^i\right).\nonumber
\end{eqnarray}
By Lemma 1,
\begin{equation}
 E\left(2\sum\limits_{i=1}^p\lambda_iQ_k^i\right)=2\sum\limits_{i=1}^{p}\lambda_i\cdot 0=0.\nonumber
\end{equation}
In view of the central limit theorem for i.i.d. random variables, (\ref{prop1:target}) is reduced to proving
\begin{equation}\label{prop1:target2}
 Var\left(2\sum\limits_{i=1}^p\lambda_iQ_k^i\right)=\sum\limits_{1\leq i,j\leq p}\lambda_i\lambda_j\Xi(i,j).
\end{equation}
By a similar argument as in Lemma 1, together with some algebraic calculations, we obtain
\begin{eqnarray*}
 &&Var\left(2\sum\limits_{i=1}^p\lambda_iQ_k^i\right)\\
 &=&4\sum\limits_{1\leq i,j\leq p}\lambda_i\lambda_jCov\left(Q_k^i,Q_k^j\right)\\
 &=&4\sum\limits_{1\leq i,j\leq p}\lambda_i\lambda_jE
 \left(\iint\limits_{S}\left(T_{\bdsm{\theta}_0}\left([a,b]\right)-Y_k\left(a,b\right)\right)
 \frac{\partial T_{\bdsm{\theta}_0}}{\partial\theta_i}\left([a,b]\right)W(a,b)\mathrm{d}a\mathrm{d}b\right)\\
 &&\left(\iint\limits_{S}\left(T_{\bdsm{\theta}_0}\left([a,b]\right)-Y_k\left(a,b\right)\right)
 \frac{\partial T_{\bdsm{\theta}_0}}{\partial\theta_j}\left([a,b]\right)W(a,b)\mathrm{d}a\mathrm{d}b\right)\\
 &=&4\sum\limits_{1\leq i,j\leq p}\lambda_i\lambda_j\iiiint\limits_{S\times S}\left\{P\left(X_1\cap[a,b]\neq\emptyset,X_1\cap[c,d]\neq\emptyset\right)
  -T_{\bdsm{\theta}_0}\left([a,b]\right)T_{\bdsm{\theta}_0}\left([c,d]\right)\right\}\\
  &&\frac{\partial T_{\bdsm{\theta}_0}}{\partial\theta_i}\left([a,b]\right)
  \frac{\partial T_{\bdsm{\theta}_0}}{\partial\theta_j}\left([c,d]\right)
  W(a,b)W(c,d)\mathrm{d}a\mathrm{d}b\mathrm{d}c\mathrm{d}d.
\end{eqnarray*}
This validates (\ref{prop1:target2}), and hence finishes the proof.\\

\end{document}